\documentclass[12pt]{amsart}
\usepackage{amssymb}
\usepackage{amsmath}
\usepackage{longtable}
\newcommand\g{{\mathfrak g}}
\newcommand\h{{\mathfrak h}}

\renewcommand{\k}{\mathfrak{k}}
\renewcommand{\b}{\mathfrak{b}}
\newcommand\gl{\mathfrak{gl}}
\renewcommand{\a}{\mathfrak{a}}

\newcommand\m{\mathfrak m}
\renewcommand\l{\mathfrak l}
\renewcommand{\L}{\mathcal{L} }

\newcommand\q{\mathfrak q}
\newcommand\z{\mathfrak z}

\renewcommand{\t}{\mathfrak{t}}

\newcommand\im{\operatorname{im}}
\newcommand\Spec{\operatorname{Spec}}

\newcommand\Pic{\operatorname{Pic}}

\newcommand\Aut{\operatorname{Aut}}

\newcommand\Supp{\operatorname{Supp}}
\newcommand\C{\mathbb C}
\newcommand\K{{\mathbb{K}}}
\newcommand\X{\mathfrak X}
\newcommand\Q{\mathbb Q}
\newcommand\R{\mathbb R}
\renewcommand\P{\mathbb P}
\newcommand\N{\mathbb N}
\newcommand\D{\mathcal D}
\newcommand\Z{\mathbb Z}

\newcommand\V{\mathcal V}

\newcommand\Rad{\operatorname{R}}
\renewcommand\sl{\mathfrak{sl}}
\newcommand\so{\mathfrak{so}}
\renewcommand\sp{\mathfrak{sp}}

\newcommand\ord{\operatorname{ord}}

\newcommand\SL{\mathop{\rm SL}\nolimits}
\newcommand\Sp{\mathop{\rm Sp}\nolimits}
\newcommand\SO{\mathop{\rm SO}\nolimits}

\newcommand\Span{\operatorname{Span}}

\newcommand{\Ad}{\mathop{\rm Ad}\nolimits}

\newcommand{\rank}{\mathop{\rm rk}\nolimits}

\renewcommand{\Ad}{\mathop{\rm Ad}\nolimits}
\newcommand\quo{/\!/}
\newtheorem{Thm}{Theorem}[section]
\newtheorem{Prop}[Thm]{Proposition}
\newtheorem{Cor}[Thm]{Corollary}
\newtheorem{Conj}[Thm]{Conjecture}
\newtheorem{Lem}[Thm]{Lemma}
\theoremstyle{definition}

\newtheorem{defi}[Thm]{Definition}
\newtheorem{Rem}[Thm]{Remark}
\newtheorem{Alg}[Thm]{Algorithm}
 \unitlength=1mm   \numberwithin{equation}{section}
\author{Ivan V. Losev}
\title{Proof of the Knop conjecture}
\begin{document}
\begin{abstract}
In this paper we prove the Knop conjecture asserting that two smooth
affine spherical varieties with the same weight monoid are
equivariantly isomorphic. We also state and prove a uniqueness
property for (not necessarily smooth) affine spherical  varieties.
\end{abstract}
\maketitle 
\section{Introduction}
Throughout the paper the base field $\K$ is algebraically closed and
of characteristic zero.

Let $G$ be a connected reductive group, $X$  an irreducible
$G$-variety. Fix a Borel subgroup $B\subset G$ and a maximal torus
$T\subset B$.

The algebra of regular functions $\K[X]$ has a natural structure of a $G$-module. 
It is known that $\K[X]$ is the sum of its finite dimensional $G$-modules. 
 By the {\it weight monoid} $\X^+_{G,X}$  of $X$ we mean the set of all highest weights of the $G$-module $\K[X]$. 
Since $\K[X]$ is an integral domain, the product of two highest vectors is non-zero, whence again
a highest vector. It follows that  $\X^+_{G,X}$ is indeed a
submonoid of the character lattice  $\X(T)$ of $T$. By  results of Popov,
[Po2],  $\X^+_{G,X}$ is finitely generated whenever $X$
is affine.

Recall  $X$ is said to be {\it spherical} iff $X$ is normal and $B$
has an open orbit in $X$. If $X$ is affine, then $X$ is spherical
iff $\K[X]$ is a multiplicity free $G$-module, that is, the
multiplicity of every irreducible module in $\K[X]$  is at most 1.
In this case $\K[X]=\bigoplus_{\lambda\in \X^+_{G,X}}V(\lambda)$ and
the monoid $\X^+_{G,X}$ is saturated, that is, $\X^+_{G,X}$ is the
intersection of a lattice with a finitely generated cone in
$\X(T)_{\otimes \Z}\Q$.

Suppose $X$ is an affine spherical variety. Let $\lambda,\mu,\nu\in
\X^+_{G,X}$ be such that $V(\nu)\subset V(\lambda)V(\mu)$ (the
product is taken in $\K[X]$). An element of the form
$\lambda+\mu-\nu$ is said to be a {\it tail} of $X$. By the {\it
tail cone} of $X$ we mean the closure of the  cone in $\t(\R)^*$
generated by all tails. Here $\t(\R)^*$ (real form of the dual space to Cartan subalgebra)
stands for $\X(T)\otimes_\Z\R$. The tail cone has a distinguished system of
generators called the system of spherical roots of $X$ and denoted
by $\Psi_{G,X}$, see Section \ref{SECTION_combinatorial} for
details.

\begin{defi}\label{defi:1}
Let $X_1,X_2$ be affine spherical varieties. We say that $X_1,X_2$
are $\X^+$-equivalent (resp., $\X^+\Psi$-equivalent) if
$\X^+_{G,X_1}=\X^+_{G,X_2}$ (resp., $\X^+_{G,X_1}=\X^+_{G,X_2},
\Psi_{G,X_1}=\Psi_{G,X_2}$).
\end{defi}

The main objective of the paper is to prove the following
assertions.
\begin{Thm}\label{Thm:1}
Any two $\X^+\Psi$-equivalent affine spherical varieties are equivariantly
isomorphic.
\end{Thm}

\begin{Thm}\label{Thm:2} Any two {\rm smooth}  $\X^+$-equivalent affine
spherical varieties are equivariantly isomorphic.
\end{Thm}

It is known that there are only finitely many systems of spherical
roots for given $G$. See \cite{Wasserman} for details.  

Theorem \ref{Thm:2} was conjectured by F.Knop. Note that it fails if the smoothness 
assumption is omitted. In fact, for $G=\SO_3$  the spherical varieties $X_0=\{(x_1,x_2,x_3)| 
x_1^2+x_2^2+x_3^2=0\}, X_1=\{(x_1,x_2,x_3)| 
x_1^2+x_2^2+x_3^2=1\}$ are $\X^+$-equivalent but not isomorphic. In this example we have $\Psi_{G,X_0}=\varnothing$, while $ \Psi_{G,X_1}$ consists of one element. 

If $G$ is a torus, i.e. $G=T$, then spherical is the same as toric.  If  $X$ is an affine toric  $T$-variety, then
there is an isomorphism of $T$-algebras $\K[X]\cong \K[\X^+_{T,X}]$.
Therefore both theorems hold. For $G$  of type $A$ (that is, when
all simple ideals of $\g$ are of type $A$)  Theorem \ref{Thm:2} was
proved  by Camus, \cite{Camus}. His approach uses Luna's
classification of spherical varieties, see~\cite{Luna5}. We do not
use his result in this paper.

We remark that for any finitely generated saturated monoid
$\X^+\subset \X(T)$ consisting of dominant weights there is an
affine spherical variety $X$ with
$\X^+_{G,X}=\X^+,\Psi_{G,X}=\varnothing$, see \cite{Popov} for
details. However, if $G$ is not a torus, then a smooth affine
spherical variety usually has at least one spherical root.

Apart of being interesting  in its own right, Theorem \ref{Thm:2} is important for the
theory of multiplicity-free Hamiltonian actions of compact groups on
compact symplectic smooth manifolds. Let us recall all necessary
definitions.

Let $K$ be a connected compact Lie group and $M$ be a  real smooth
manifold equipped with a symplectic form $\omega$. An action $K:M$
is called Hamiltonian if it preserves $\omega$ and is equipped with
a {\it moment map} $\mu$, that is, a $K$-equivariant map
$\mu:M\rightarrow \k^*$ satisfying
$$\omega(\xi_*x,v)=\langle d_x\mu(v),\xi\rangle, \forall x\in M,\xi\in\k,v\in T_xM.$$
Here  $\xi_*x$ denotes the  vector field at $x$ corresponding to
$\xi$. By  a Hamiltonian $K$-manifold we mean a symplectic manifold
equipped with a Hamiltonian action of $K$.

Choosing a $K$-invariant scalar product on $\k$, one identifies
$\k^*$ with $\k$ and considers $\mu$ as a map $M\rightarrow\k$. Fix
a Cartan subalgebra $\t\subset\k$ and a positive Weyl chamber
$\t_+\subset\t$. Define the map $\psi:M\rightarrow \t_+$ by
$\psi(x)=K\mu(x)\cap\t_+$.

Recall that a compact  Hamiltonian $K$-manifold $M$ is called {\it
multiplicity-free} if it satisfies the following equivalent
conditions:
\begin{enumerate}
\item A general $K$-orbit in $M$ is a coisotropic submanifold.
\item Any fiber of $\psi$ is a single $K$-orbit.
\item The algebra $C^\infty(M)^K$ is commutative with respect to the
 Poisson bracket induced from $C^\infty(M)$.
\end{enumerate}
The proof of equivalence is similar, for example, to the proof of
Proposition A.1 in \cite{Woodward}.

An important and interesting problem is to classify (in reasonable
terms) all compact multiplicity-free Hamiltonian $K$-manifolds. To
solve this problem we need to introduce certain simple invariants
associated with such a manifold.  Firstly, by the Kirwan theorem,
\cite{Kirwan}, $\psi(M)$ is a convex polytope called the {\it moment
polytope} of $M$.  This is the first invariant we need. The second
one is the stabilizer $K_x$ for $x\in \mu^{-1}(\eta)$, where $\eta$
is a general element of $\psi(M)$. This stabilizer does not depend on choices of 
$\eta$ and $x$  so
is determined uniquely. It is called   the {\it
principal isotropy subgroup} of $M$. Note that for an arbitrary action one can define a principal 
isotropy group only up to conjugacy but in our situation we have a "distinguished" general point $x$, that lying in
$\mu^{-1}(C)$. 

\begin{Conj}[Delzant]
Any two multiplicity-free compact Hamiltonian $K$-manifolds $M_1,M_2$ with
the same moment polytope and principal isotropy subgroup are
isomorphic (that is, there is a $K$-equivariant symplectomorphism
$\varphi:M_1\rightarrow M_2$ commuting with the moment maps).
\end{Conj}

Currently, the conjecture is proved only in some special cases: see
\cite{Delzant1},\cite{Delzant2},\cite{Woodward}.  F. Knop reduced
the Delzant conjecture to Theorem \ref{Thm:2} (unpublished). So  this paper
completes the proof of the Delzant conjecture (modulo Knop's reduction).

Theorem \ref{Thm:1} also has some applications to symplectic
geometry. Using this theorem, one proves a certain uniqueness result
for invariant K\"{a}hler structures on a given multiplicity free
compact Hamiltonian $K$-manifold, Theorem \ref{Thm:8.1}.

Let us describe briefly the content of the paper.  Section
\ref{SECTION_Notation} contains conventions used in the paper and
the list of notation. In Section \ref{SECTION_combinatorial} we
recall the definitions and some properties of important combinatorial
invariants of spherical varieties. These invariants are the Cartan
space $\a_{G,X}$, the weight lattice $\X_{G,X}$, the valuation cone
$\V_{G,X}$, the system of spherical roots $\Psi_{G,X}$, and the set
of $B$-divisors $\D_{G,X}$ equipped with certain two maps.

Section \ref{SECTION_monoid} contains some auxiliary results
concerning affine spherical varieties and weight monoids. Further,
we  state there two auxiliary statements~-- Theorems \ref{Thm:5.1},
\ref{Thm:6.1}. Theorem \ref{Thm:5.1} asserts, roughly speaking, that
affine spherical varieties $X_1,X_2$ have the same set of
$B$-divisors provided they are $\X^+\Psi$-equivalent or smooth and
$\X^+$-equivalent. Theorem \ref{Thm:6.1} states that
$\Psi_{G,X_1}=\Psi_{G,X_2}$ provided $X_1,X_2$ are smooth and
$\X^+$-equivalent. In the end of Section \ref{SECTION_monoid} we
deduce Theorems \ref{Thm:1},\ref{Thm:2} from Theorems
\ref{Thm:5.1},\ref{Thm:6.1} and results of \cite{unique}.

Section \ref{SECTION_Reductions} is devoted to reduction procedures,
which are based on the local structure theorem and play a crucial
role  in the proofs of Theorems \ref{Thm:5.1},\ref{Thm:6.1}. These
proofs are presented in Sections
\ref{SECTION_divisors},\ref{SECTION_roots}. At the end of Section
\ref{SECTION_divisors} we also give an algorithm  recovering
$\D_{G,X}$ from $\X^+_{G,X}$ and $\Psi_{G,X}$.

Finally in Section \ref{SECTION_Kahler} we  prove a uniqueness
result for invariant K\"{a}hler structures on a given compact
multiplicity free Hamiltonian manifold.

{\bf Acknowledgements}. This paper was partially written during
author's stay in Rutgers University, New Brunswick,  in the
beginning of 2007. I would like to thank this institution and
especially  Professor F. Knop for hospitality. I am also grateful to
B. Van Steirteghem for some useful remarks on an earlier version of
this text.

\section{Notation, terminology and
conventions}\label{SECTION_Notation} If an algebraic group is
denoted by a capital Latin letter, then we denote its Lie algebra by
the corresponding small German letter.

Recall that we fix a Borel subgroup $B\subset G$ and a maximal torus
$T\subset B$. This allows us to define the root system $\Delta(\g)$,
the Weyl group $W(\g)$, and the system of simple roots $\Pi(\g)$ of
$\g$. For simple roots and fundamental weights of $\g$ we use the
notation of \cite{VO}. Put $U=\Rad_u(B)$. 

Note that the character groups $\X(T),\X(B)$
are naturally identified. We recall that the character group $\X(H)$ of an algebraic
group $H$, by definition, consists of all algebraic group homomorphisms $H\rightarrow \K^
\times$, where $\K^\times$ is the one-dimensional torus. If $H$ is connected, then $\X(H)$
is a lattice, so we call it the {\it character lattice}.

Since $G$ is reductive, there is a $G$-invariant symmetric form
$(\cdot,\cdot)$, whose restriction to $\t(\Q)$ is positive
definite. Note that $(\cdot,\cdot)$ is nondegenerate on the Lie
algebra of any reductive subgroup of $G$. We fix such a form and
identify $\g$ with $\g^*$, $\t$ with $\t^*$.

If $X_1,X_2$ are  $G$-varieties, then we write  $X_1\cong^G X_2$
when $X_1,X_2$ are $G$-equivariantly isomorphic.  Note that if
$V_1,V_2$ are $G$-modules, then $V_1\cong^G V_2$ iff $V_1,V_2$ are
isomorphic as $G$-modules.

Let $Q$ be a parabolic subgroup of $G$ containing either $B$ or
$B^-$. There is a unique Levi subgroup of $Q$ containing $T$, we
call it the {\it standard} Levi subgroup of $Q$. For $\Sigma\subset
\Pi(\g)$ we denote by $P_\Sigma$ the  parabolic subgroup of $G$
whose Lie algebra is generated by $\b$ and the root subspaces
corresponding to $-\alpha$ with $\alpha\in\Sigma$.

\begin{longtable}{p{5.5cm} p{10.3cm}}
$A^{(B)}_\lambda$&$=\{a\in A| b.a=\lambda(b)a, \forall b\in B\}$.
\\$A^{(B)}$&$=\cup_{\lambda\in \X(B)}A^{(B)}_\lambda$.
\\
$\a_{G,X}$& the Cartan space of a spherical $G$-variety $X$, see Section \ref{SECTION_combinatorial}
for the definition.
\\ $\D_{G,X}$& the set of  $B$-divisors of a spherical
$G$-variety, see Section \ref{SECTION_combinatorial}
for the definition.
\\ $\D_{G,X}(\alpha)$& $=\{D\in\D_{G,X}| P_\alpha\not\subset G_D\}$.
\\ $\D_{G,X}^G$&$=\{D\in \D_{G,X}| G_D=G\}$.
\\ $f_\lambda$& a nonzero element of $\K(X)^{(B)}_\lambda$.
\\ $(f)$& the zero divisor of a rational function $f$.
\\ $(G,G)$& the derived subgroup of a group
$G$\\
$[\g,\g]$& the derived subalgebra of a Lie algebra $\g$.
\\ $G*_HV$& the homogeneous bundle over $G/H$ with a fiber $V$.\\ $G_y$& the stabilizer of  $y$ under an
action of $G$.
\\ $\K^\times$& the one-dimensional torus.
\\ $N_G(H)$& $=\{g\in G| gHg^{-1}=H\}$.
\\ $\Pic(X)$& the Picard group of a variety $X$.
\\ $\Rad_u(G)$& the unipotent radical of an algebraic group $G$.
\\ $\rank_G(X)$& $=\rank\X_{G,X}$.
\\ $\Span_{A}(M)$&$=\{a_1m_1+\ldots+a_km_k, a_i\in A, m_i\in M\}$.
\\ $\Supp(\gamma)$& the support of  $\gamma\in\Span_\Q(\Pi(\g))$, that is, the set $\{\alpha\in \Pi(\g)|n_\alpha\neq 0\}$, where
$\gamma=\sum_{\alpha\in\Pi(\g)}n_\alpha\alpha$.
\\ $V(\mu)$& the irreducible module with  highest weight $\mu$.\\
$X^G$&$=\{x\in X| gx=x, \forall g\in G\}$.\\
$\X(G)$& the character group of an algebraic group $G$.\\
$\X_{G,X}$& the weight lattice of a spherical $G$-variety $X$, see Section \ref{SECTION_combinatorial}
for definition.
\\  $\X^+_{G,X}$&
the weight monoid of a spherical $G$-variety $X$.\\
 $\#X$& the cardinality of a
set $X$.
\\ $W(\g)$& the Weyl group of a reductive Lie algebra $\g$.
\\   $Z_G(\h)$&=$\{g\in G| \Ad(g)|_\h=id\}$\\  $\z_\g(\h)$&$=\{\xi\in\g| [\xi,\h]=0\}$.
\\ $\alpha^\vee$& the dual root corresponding to a root $\alpha$.
\\$\langle\alpha,v\rangle$& the pairing of  elements $\alpha,v$ of
dual vector spaces.
\\  $\Delta(\g)$& the root system of a reductive Lie algebra $\g$.
\\ $\Pi(\g)$& the system of simple roots of $\g$ associated with
$B$.
\\ $\Pi(\g)_X^a,\ldots,\Pi(\g)_X^d$& the set of all simple roots of
type a),\ldots,d) for a spherical $G$-variety $X$, see Section \ref{SECTION_combinatorial}
for the definition.
\\ $\Psi_{G,X}$& the system of spherical roots of a spherical
$G$-variety $X$, see Section \ref{SECTION_combinatorial}
for the definition.
\\ $\varphi_D$& the vector in $\a_{G,X}^*$ associated with $D\in
\D_{G,X}$, see Section \ref{SECTION_combinatorial}
for the definition.
\end{longtable}

\section{Combinatorial invariants of spherical
varieties}\label{SECTION_combinatorial} In this section $G$ is a
connected reductive algebraic group and
 $X$ is a spherical $G$-variety.

By the {\it weight lattice} of $X$ we mean the set
$\X_{G,X}:=\{\mu\in\X(T)| \K(X)^{(B)}_\mu\neq \{0\}\}$. This is a
sublattice in $\X(T)$. 
Note that $\dim\K(X)^{(B)}_\mu=1$
for any $\mu\in \X_{G,X}$.

Define the {\it Cartan space} of $X$ by
$\a_{G,X}:=\X_{G,X}\otimes_\Z\Q$. This is a subspace in $\t(\Q)^*$.
Note  that $\a_{G,X}$ is equipped with a positive definite
symmetric bilinear form (the restriction of the form on $\t^*(\Q)$).

Let $v$ be a $\Q$-valued discrete valuation of $\K(X)$. One defines
the element $\varphi_v\in \a_{G,X}^*$ by
\begin{equation*}
\langle\varphi_v,\mu\rangle=v(f_\mu), \forall \mu\in\X_{G,X}.
\end{equation*}

  It is known, see \cite{Knop5}, that the restriction of the map $v\mapsto \varphi_v$
to the set of all $G$-invariant $\Q$-valued discrete valuations is
injective. Its image is a finitely generated cone in $\a_{G,X}^*$.
We denote this image by $\V_{G,X}$ and call it the {\it valuation
cone}.

Let $X$ be, in addition, affine and $\mathcal{T}$ denote the tail
cone of $X$. It is known, see, for example, \cite{Knop5}, Lemma 5.1,
that  $-\V_{G,X}$ is  the dual cone of $\mathcal{T}$. In  turn
$\mathcal{T}$ is the dual cone of $-\V_{G,X}$.

 Moreover, the valuation cone is
a Weyl chamber for a (uniquely determined) group $W_{G,X}$ generated
by reflections, see \cite{Knop3}, Theorem 7.4. This group is called
the {\it Weyl group} of $X$.

Denote by $\Psi_{G,X}$  the set of primitive elements $\alpha\in
\X_{G,X}$  such that $\ker\alpha\subset \a_{G,X}^*$ is a wall of
$\V_{G,X}$ and $\alpha$ is nonpositive on $\V_{G,X}$. It is clear
from construction that $\Psi_{G,X}$ is a system of simple roots in a
certain root system with Weyl group $W_{G,X}$. An element of
$\Psi_{G,X}$ is called a {\it spherical root} of $X$. So if $X$ is
affine, then its tail cone is generated by $\Psi_{G,X}$.

By a $B$-divisor we mean a prime $B$-stable divisor.  Let $\D_{G,X}$
denote the set of all  $B$-divisors on $X$. We write $\varphi_D$
instead of $\varphi_{\ord_D}.$ 
Further, set $G_D:=\{g\in G| gD=D\}$. Clearly, $G_D$ is a parabolic
subgroup of $G$ containing $B$. For a subset $\D\subset\D_{G,X}$ we
put $G_\D:=\cap_{D\in\D}G_\D$.

Let $\alpha\in \Pi(\g)$. Set $\D_{G,X}(\alpha):=\{D\in \D_{G,X}|
P_\alpha\not\subset G_D\}, \D_{G,X}^G:=\{D\in \D_{G,X}| G_D=G\}$.
Clearly,
$\D_{G,X}=\D_{G,X}^G\sqcup\bigcup_{\alpha\in\Pi(\g)}\D_{G,X}(\alpha)$.

Now we are going to recall Luna's results (\cite{Luna4}, see also
\cite{Luna5}, Section 2) concerning the structure of the sets
$\D_{G,X}(\alpha)$ and the vectors $\varphi_D,D\in
\D_{G,X}(\alpha)$.

\begin{Prop}[\cite{Luna4}, Proposition 3.4, \cite{Luna5}, Lemma 6.4.2]\label{Prop:2.3.1}  For $\alpha\in \Pi(\g)$ exactly one of the following possibilities
takes place:
\begin{itemize}
\item[(a)] $\D_{G,X}(\alpha)=\varnothing$.
\item[(b)] $\alpha\in \Psi_{G,X}$. Here
$\D_{G,X}(\alpha)=\{D^+,D^-\}$,
$\varphi_{D^+}+\varphi_{D^-}=\alpha^\vee|_{\a_{G,X}},
\langle\varphi_{D^\pm},\alpha\rangle=1$.
\item[(c)]  $2\alpha\in \Psi_{G,X}$. In this case $\D_{G,X}(\alpha)=\{D\}$ and $\varphi_{D}=\frac{1}{2}\alpha^\vee|_{\a_{G,X}}$.
\item[(d)] $\Q\alpha\cap \Psi_{G,X}=\varnothing, \D_{G,X}(\alpha)=\{D\}$ and $\varphi_{D}=\alpha^\vee|_{\a_{G,X}}$.
\end{itemize}
\end{Prop}

We say that a root $\alpha\in \Pi(\g)$ is {\it of type} a) (or
b),c),d)) if the corresponding possibility takes place for $\alpha$.
The set of all simple roots of type a),\ldots, d) is denoted by
$\Pi(\g)_X^a,\ldots,\Pi(\g)_X^d$.

\begin{Prop}[\cite{Luna5}, Proposition 3.2]\label{Prop:2.3.2}
Let $\alpha,\beta\in \Pi(\g)$. If
$\D_{G,X}(\alpha)\cap\D_{G,X}(\beta)\neq\varnothing$, then exactly
one of the following possibilities takes place:
\begin{enumerate}
\item $\alpha,\beta\in\Pi(\g)_X^b$ and
$\#\D_{G,X}(\alpha)\cap\D_{G,X}(\beta)=1$.
\item $\alpha,\beta\in \Pi(\g)_X^d, \langle\alpha^\vee,\beta\rangle=0,\alpha^\vee-\beta^\vee|_{\a_{G,X}}=0$
and $\alpha+\beta=\gamma$ or $2\gamma$ for some
$\gamma\in\Psi_{G,X}$.
\end{enumerate}
Conversely, if $\alpha,\beta\in\Pi(\g)$ are as in (2), then
$\D_{G,X}(\alpha)=\D_{G,X}(\beta)$.
\end{Prop}

The following lemma is proved analogously to Lemma 4.1.12 from
\cite{unique}.

\begin{Lem}\label{Lem:2.3.2}
Let $X$ be a spherical $G$-variety, and $\alpha\in \Psi_{G,X}$ have
one of the following forms:
\begin{enumerate}
\item $\alpha=\alpha_1$, where $\alpha_1\in \Pi(\g)$.
\item $\alpha=2\alpha_1$, where $\alpha_1\in\Pi(\g)$.
\item $\alpha=k(\alpha_1+\alpha_2)$, where $\alpha_1,\alpha_2$ are
orthogonal elements of $\Pi(\g)$ and $k=1$ or $\frac{1}{2}$.
\end{enumerate}
If $D\in \D_{G,X}\setminus\D_{G,X}(\alpha_1)$, then
$\langle\varphi_D,\alpha\rangle\leqslant 0$.
\end{Lem}

\begin{Cor}\label{Cor:2.3.8}
Let $\alpha\in \Pi(\g), D\in \D_{G,X}(\alpha)$.
\begin{enumerate}
\item If $\alpha\in \Pi(\g)_X^b\cup\Pi(\g)_X^c$, then $G_D=P_{\Pi(\g)\setminus A}$,
 where $A:=\{\beta\in \Pi(\g)|\langle\beta,\varphi_D\rangle=1\}$. In
 particular, if $\alpha\in \Pi(\g)_X^c$, then $A=\{\alpha\}$.
\item If $\alpha\in \Pi(\g)_X^d$, then either $G_D=P_{\Pi(\g)\setminus
\{\alpha\}}$ or $G_D=P_{\Pi(\g)\setminus \{\alpha,\beta\}}$, where
$\beta$ is the unique simple root such that
$\langle\alpha^\vee,\beta\rangle=0,
\alpha^\vee-\beta^\vee|_{\a_{G,X}}=0$ and $\alpha+\beta$ is a
positive multiple of a spherical root.
\end{enumerate}
\end{Cor}
\begin{proof}
The first assertion follows from Proposition~\ref{Prop:2.3.2} and
Lemma \ref{Lem:2.3.2}. In assertion 2 one only needs to prove that
$\beta$ is unique. Assume the contrary, let $\beta,\gamma$ be such
that $\D_{G,X}(\alpha)=\D_{G,X}(\beta)=\D_{G,X}(\gamma)$. Note that
$k_1(\alpha+\beta), k_2(\alpha+\gamma), k_3(\beta+\gamma)\in
\Psi_{G,X}$ for some $k_1,k_2,k_3>0$. But
$\langle\alpha^\vee-\beta^\vee, \beta+\gamma\rangle\neq 0$, which
contradicts Proposition \ref{Prop:2.3.2}.
\end{proof}

\section{Some remarks on affine spherical varieties and weight
monoids}\label{SECTION_monoid} Throughout this section $X$ is an
affine spherical variety.

The definition of the weight monoid $\X^+_{G,X}$ given in 
the introduction can be rewritten as $\X^+_{G,X}=\{\lambda\in \X(T)| \K[X]^{(B)}_\lambda\neq \{0\}\}$.
The following lemma is easy. A proof can be found, for instance, in \cite{unique}, Lemma 3.6.2.

\begin{Lem}\label{Lem:2.2.5}
$\X_{G,X}=\Span_\Z(\X_{G,X}^+)$.
\end{Lem}

\begin{Lem}\label{Cor:2.2.8}
If $X_1,X_2$ are $\X^+$-equivalent affine spherical $G$-varieties,
then $\dim X_1=\dim X_2$.
\end{Lem}
\begin{proof}
Since $\K[X_1]^U\cong\K[\X^+_{G,X_i}]\cong \K[X_2]^U$, the claim
follows easily from results of \cite{Popov}.
\end{proof}

The following proposition follows from the Luna slice theorem, see \cite{KvS}, Corollary 2.2,
for details. 

\begin{Prop}
If $X$ is smooth, then $X\cong^G G*_HV$, where $H$ is a reductive
subgroup of $G$ and $V$ is an $H$-module.
\end{Prop}

\begin{Rem}\label{Rem:1.4.1}
Note that  $G$ acts on the set of all pairs $(H,V)$:
$g.(H,V)=(H',V')$, where  $H'=gHg^{-1}$ and there is a linear
isomorphism $\iota:V\rightarrow V'$ such that
$(ghg^{-1})\iota(v)=\iota(gv)$ for all $h\in H$. One easily shows,
see, for instance, \cite{unique}, Lemma 3.6.6, that $G*_{H}V\cong^G
G*_{H'}V'$ iff $(H,V)\sim_G
(H',V')$.  
So if $V_0$ is a $G$-module and
$V_0\times (G*_HV)\cong V_0\times (G*_{H'}V')$, then $G*_HV\cong
G*_{H'}V'$.
\end{Rem}

The following lemma follows directly from highest weight theory.

\begin{Lem}\label{Lem:1.1.1}
  A simple normal subgroup
$G_1\subset G$ acts trivially on $X$ iff  $\langle\X^+_{G,X}, T\cap
G_1\rangle=1$. An element $t\in Z(G)$ acts trivially on $X$ iff
$\langle\X^+_{G,X},t\rangle=1$.
\end{Lem}

\begin{defi}\label{Def:1.1.8}
A  $G$-variety $X$ is said to be {\it decomposable} if there exist
nontrivial connected normal subgroups $G_1,G_2\subset G$ and
$G_i$-varieties $X_i, i=1,2,$ such that $G$ is decomposed into a
locally direct product of $G_1,G_2$ (that is, $G=G_1G_2$ and $G_1\cap G_2$
is finite) and $X\cong^{G_1\times G_2}
X_1\times X_2$. Under these assumptions we say that the pair
$G_1,G_2$ {\it decomposes} $X$.
\end{defi}

\begin{Lem}\label{Lem:1.1.7}
Let  $G=G_1G_2$ be a decomposition into a locally direct product.
Then the following conditions are equivalent.
\begin{enumerate}
\item The pair $(G_1,G_2)$ decomposes $X$.
\item $\X^+_{G_1\times G_2,X}=\Gamma_1+\Gamma_2$, where $\Gamma_i\subset \X(T\cap G_i),
i=1,2$.
\end{enumerate}
In particular, if affine spherical varieties $X_1,X_2$ are
$\X^+$-equivalent, and the pair $(G_1,G_2)$ decomposes $X_1$, then
it decomposes $X_2$.
\end{Lem}
\begin{proof}
Essentially, this lemma was proved in the proof of Lemma 3.6.4 in
\cite{unique}. In order to make the present paper more self-contained 
we present an argument below. 

Clearly, $(1)\Rightarrow (2)$. Let us check the opposite implication. 
Set $X^i:=X\quo G_i$.
From 
highest weight theory one easily deduces that
$\X_{G_i,X^i}^+=\X_{G,X}^+\cap \X(T\cap G^i),i=1,2$. Thus
$\X^+_{G,X}=\X^+_{G_1,X^1}+\X^+_{G_2,X^2}$. In other words,
$\K[X]^{U}=\K[X^1]^{U\cap G_1}\otimes
\K[X^2]^{U\cap G_2}$. It follows from  highest weight
theory that $\K[X]=\K[X^1]\otimes \K[X^2]$.
\end{proof}

To prove Theorems \ref{Thm:1},\ref{Thm:2} we need the following two
theorems.

\begin{Thm}\label{Thm:5.1}
Let $X_1,X_2$ be  affine spherical varieties. If $X_1,X_2$ are
$\X^+\Psi$-equivalent or $X_1,X_2$ are $\X^+$-equivalent and smooth,
then  there is a bijection $\iota:\D_{G,X_1}\rightarrow \D_{G,X_2}$
such that $\varphi_{\iota(D)}=\varphi_D$, $G_{\iota(D)}=G_D$.
\end{Thm}

\begin{Thm}\label{Thm:6.1}
Let $X_1,X_2$ be smooth $\X^+$-equivalent affine spherical
varieties. Then $\Psi_{G,X_1}=\Psi_{G,X_2}$.
\end{Thm}

\begin{proof}[Proof of Theorems \ref{Thm:1},\ref{Thm:2} modulo Theorems \ref{Thm:5.1},\ref{Thm:6.1}]
Let $X_1,X_2$ be either $\X^+\Psi$-equivalent or smooth and
$\X^+$-equivalent affine spherical varieties.  Let $X^0_i,i=1,2,$
denote the open $G$-orbit of $X_i$. Clearly,
$\X_{G,X^0_i}=\X_{G,X_i},\V_{G,X^0_i}=\V_{G,X_i},\D_{G,X^0_i}=\D_{G,X_i}\setminus
\D_{G,X_i}^G, i=1,2$. Thanks to Lemma \ref{Lem:2.2.5}, Theorems
\ref{Thm:5.1}, \ref{Thm:6.1}, one can apply Theorem 1 from
\cite{unique} to $X_1^0,X_2^0$ and obtain that $X_1^0\cong^G X_2^0$.
It follows from \cite{unique}, Proposition 3.6.5, that $X_1\cong^G
X_2$.
\end{proof}

In the proof of Theorems \ref{Thm:5.1},\ref{Thm:6.1} we may (and
will) assume that Theorems \ref{Thm:1},\ref{Thm:2} are already
proved for all groups $G'$ and spherical ($\X^+\Psi$-equivalent or
smooth and $\X^+$-equivalent) $G'$-varieties $X_1',X_2'$ provided
one of the following assumptions holds:
\begin{itemize}
\item[(A1)] $\dim G'<\dim G$.
\item[(A2)] $G=G'$ and $\min_{i=1,2}(\#\D_{G,X_i}^G)>\min_{i=1,2}(\#\D_{G,X_i'}^G)$.
\end{itemize}

\section{Reductions}\label{SECTION_Reductions}
Our reductions are based on the local structure theorem.   First
variants of this theorem were proved independently in
\cite{BLV},\cite{Gr}. We give here a version due to Knop.

\begin{Thm}[\cite{Knop3}, Theorem 2.3 and Lemma 2.1]\label{Thm:2.1.1}
Let $X$ be an irreducible normal $G$-variety, $D$ an effective
$B$-stable Cartier divisor on $X$, $L_D$ the standard Levi subgroup
of $G_D$. Then there exists a closed $L_D$-subvariety $\Sigma\subset
X^0:=X\setminus D$ such that the morphism
$\Rad_u(G_D)\times\Sigma\rightarrow X^0, [g,s]\mapsto gs,$ is a
$G_D$-equivariant isomorphism. Here $\Rad_u(G_D)$ acts by left
translations on itself and trivially on $\Sigma$, while $L_D$ acts
by conjugations on $\Rad_u(G_D)$ and initially on $\Sigma$.
\end{Thm}

A subvariety $\Sigma$ satisfying the claim of the previous theorem
is said to be  a {\it section} of $X$ associated with $D$.

\begin{Rem}\label{Rem:2.1.2}
Being the quotient for the action $\Rad_u(G_D):X^0$, the
$L_D$-variety $\Sigma$ depends  only on the support of $D$.
 If
$D=\sum_{i=1}^k a_i D_i, a_i\in \N$, we denote the $L_D$-variety
$\Sigma$ by $X(D_1,\ldots, D_k)$ or $X(\{D_1,\ldots, D_k\})$. Note
also that $\Sigma$ is smooth (spherical) provided $X$ is.
\end{Rem}

Till the end of this subsection $X$ is an affine spherical variety.

Here is the first version of our reduction procedure.

\begin{Prop}\label{Prop:2.2.1}
Choose a subset $\D\subset \D_{G,X}$ and let $M$ be the standard
Levi subgroup of $G_\D$. Suppose there is a Cartier divisor of the
form $\sum_{D\in \D}a_D D, a_D>0,$ (which is the case, for instance,
when $X$ is smooth). Then
\begin{enumerate}
\item $X(\D)$ is an affine
$M$-variety.
\item  The map $\iota:\D_{M,X(\D)}\rightarrow
\D_{G,X}, D\mapsto \Rad_u(G_\D)\times D,$ is an injection with image
$\D_{G,X}\setminus \D$. Furthermore,
$\varphi_{D}=\varphi_{\iota(D)}, G_{\iota(D)}\cap M=M_{D}$ for any
$D\in \D_{M,X(\D)}$.
\item $\X_{M,X(\D)}=\X_{G,X},
\Psi_{M,X(\D)}=\Psi_{G,X}\cap\Span_\Q(\Delta(\m))$.
\item $\X^+_{M,X(\D)}=\{\chi\in \X_{G,X}|
\langle\varphi_D,\chi\rangle\geqslant 0, \forall
D\in\D_{G,X}\setminus \D\}.$
\end{enumerate}
\end{Prop}
\begin{proof}
To prove assertion 1 note that, being a complement to a Cartier divisor in an affine
variety, $X^0$ is affine. Being isomorphic to a closed
subvariety of $X^0$, the variety $\Sigma$ is affine. Assertions 2,3
follow  from \cite{unique}, Lemma  3.5.5. Assertion 4
 follows from assertion 2 and \cite{unique}, Lemma 3.6.2.
\end{proof}

In   the second version of our reduction procedure we do not need to
know $\varphi_D,D\not\in\D, G_\D$.

\begin{Prop}\label{Prop:2.2.2}
In  the above notation  suppose $(f_\mu)=\sum_{D\in\D}a_D D$ for
some positive $a_D$. Then
$G_\D=BG_\mu,\X^+_{G_\mu,X(\D)}=\X^+_{G,X}+\Z\mu$ and the image of
$\D_{G_\mu,X(\D)}$ in $\D_{G,X}$ (see assertion 2 of
Proposition~\ref{Prop:2.2.1}) coincides with
$\{D\in\D_{G,X}|\langle\mu,\varphi_D\rangle=0\}$.
\end{Prop}
\begin{proof}
An element $g\in \K(X)$ lies in $\K[X\setminus (f)]$ iff
$g=\frac{g_1}{f^n}$ for some $n\in\N, g_1\in K[X]$. Clearly,
$g\in\K[X^0]^{(B)}$ iff $g_1\in \K[X]^{(B)}$.  The equality for
$\X^+_{G_\mu,X(\D)}$ stems now from the natural isomorphism
$\K[X^0]^U\cong \K[X(\D)]^{U\cap G_\mu}$. The description of
$\D_{G_\mu,X(\D)}$ is clear.
\end{proof}

Under the assumptions of the previous proposition, we put
$X(\mu):=X(\D)$. Note that $X(\mu)\neq X$ iff $f_\mu$ is not
invertible in $\K[X]$ iff $\mu\in \X^+_{G,X}\setminus -\X^+_{G,X}$.

\begin{Cor}\label{Cor:2.2.4}
Let  $\mu\in \X^+_{G,X}$ be such that $\X^+_{G,X}+\Z\mu=\X_{G,X}$
(such $\mu$ always exists). Then $\Pi(\g)_X^a=\{\alpha\in \Pi(\g)|
\langle\alpha^\vee,\mu\rangle=0\}$.
\end{Cor}
\begin{proof}
Any $f\in \K[X(\mu)]^{(B\cap G_\mu)}$ is invertible. Thus
$\D_{G(\mu),X(\mu)}=\varnothing$. From Proposition \ref{Prop:2.2.2}
it follows that $G_{\D_{G,X}}=BG_\mu$. But $
G_{\D_{G,X}}=P_{\Pi(\g)^a_X}$.
\end{proof}

The following statement stems directly from Propositions
\ref{Prop:2.2.1},\ref{Prop:2.2.2}.

\begin{Cor}\label{Cor:2.2.3}
Let $X_1,X_2$ be affine spherical $G$-varieties.
\begin{enumerate}
\item If $X_1,X_2$ are $\X^+$-equivalent (resp., $\X^+\Psi$-equivalent), then for any $\mu\in
\X^+_{G,X}\setminus -\X^+_{G,X}$ the $G_\mu$-varieties
$X_1(\mu),X_2(\mu)$ are $\X^+$-equivalent (resp.,
$\X^+\Psi$-equivalent).
\item Suppose $X_1,X_2$ are smooth and $\X^+$-equivalent, $D_1\in \D_{G,X_1}, D_2\in\D_{G,X_2}$ are such that
$G_{D_1}=G_{D_2}$ and  $\{\varphi_{D},D\in
\D_{G,X_1}\setminus\{D_1\}\}=\{\varphi_{D},D\in
\D_{G,X_2}\setminus\{D_2\}\}$. Then $X_1(D_1), X_2(D_2)$ are
$\X^+$-equivalent (as $M$-varieties, where $M$ is the standard Levi
subgroup in $G_D$).
\end{enumerate}
\end{Cor}

\begin{Rem}\label{Rem:2.2.12}
Let $X,Y$ be smooth spherical varieties and $\varphi:X\rightarrow Y$
be a smooth surjective morphism with irreducible fibers. Then
$\varphi$ induces an embedding $\varphi^*:\D_{G,Y}\hookrightarrow
\D_{G,X},D\mapsto\varphi^{-1}(D)$. Note that $G_D=G_{\varphi^*(D)}$.
\end{Rem}

In general, it is difficult to understand the structure of $X(\D)$
as a homogeneous vector bundle. However, there is a special case
when the description is easy.

\begin{Lem}\label{Lem:2.2.7}
Let $Q^-$ be a parabolic subgroup of $G$ containing $B^-$, $M$ the
standard Levi subgroup of $G$.  Suppose $X=G*_HV$, where $H\subset
M$. Denote by $\pi$ the natural homomorphism $X=G*_HV\rightarrow
G/Q^-, [g,v]\mapsto gQ^-,$ and set $\D:=\pi^*(\D_{G,G/Q^-})$. Then
$G_{\D}=BM$ and $X(\D)\cong^M \Rad_u(\q^-)\times M*_HV$.
\end{Lem}
\begin{proof}
Note that $X(\D)\cong^M \pi^{-1}(eQ^-)\cong^M Q^-*_HV\cong^M
Q^-*_M(M*_HV)\cong^M \Rad_u(Q^-)\times M*_HV$. The exponential mapping defines an 
$M$-equivariant isomorphism $\Rad_u(\q^-)\rightarrow \Rad_u(Q^-)$. 
\end{proof}

\section{Correspondence between
$B$-divisors}\label{SECTION_divisors} The goal of this section is to
prove Theorem \ref{Thm:5.1}. We  assume that $G$ is not a torus.
Throughout the section $X$ is an affine spherical variety and the
action $G:X$ is assumed to be locally effective, that is, its kernel is
finite. We write
$\a,\X^+,\X,\Psi,\D, \Pi^a,\ldots,\Pi^d$ instead of
$\a_{G,X},\X^+_{G,X},\X_{G,X},\Psi_{G,X},\D_{G,X},$
$\Pi(\g)_X^a,\ldots, \Pi(\g)_X^d$.

Let us briefly describe  the scheme of proof. The key idea in the
proof is to use Proposition \ref{Prop:2.2.2}. This proposition
allows to recover (partially) elements of $\D_{G,X_i}$ satisfying
$\langle\varphi_D,\mu\rangle=0$ for some $\mu\in
\X^+_{G,X_i}\setminus -\X^+_{G,X_i}$. This motivates us to define a
certain subset of {\it hidden} elements of $\D_{G,X_i}$ (Definition
\ref{Def:2.3.3}), those that can not be recovered by using
Proposition \ref{Prop:2.2.2}. Then we study some properties of
hidden divisors, Lemma \ref{Lem:2.3.7}, Propositions
\ref{Prop:3.1.4},\ref{Prop:2.4.1}. The former proposition deals with
the case when all divisors in $\D_{G,X_i}(\alpha)$ are hidden for
some $i$. Proposition \ref{Prop:2.4.1} describes the set of hidden
divisors for smooth $X_i$. In the proof of Theorem \ref{Thm:5.1} we
first construct a certain bijection between the sets of nonhidden
divisors. Its existence  is deduced essentially from Proposition
\ref{Prop:2.2.2}. Then we show (and this is the most complicated and
technical part of the proof) that this bijection can be extended in
the required way to the whole sets of divisors.

\begin{Lem}\label{Lem:5.2}
Suppose $X$ is a spherical $G$-variety and the action $G:X$ is
locally effective. If $\Pi^a\cup\bigcup_{\alpha\in
\Psi}\Supp(\alpha)=\Pi(\g)$, then
$\bigcup_{\alpha\in\Psi}\Supp(\alpha)=\Pi(\g)$.
\end{Lem}
\begin{proof}
Assume the contrary. It follows from Corollary \ref{Cor:2.2.4} that
$\Pi(\g_1)\not\subset \Pi^a$ for any simple ideal $\g_1\subset \g$.
Therefore there is $\beta\in \Pi^a$ such that $\beta\not\in
\Supp(\alpha)$ for any $\alpha\in \Psi$ but $\beta$ is adjacent to
$\Supp(\alpha)$ for some $\alpha\in \Psi$. It follows that
$\langle\beta^\vee,\alpha\rangle\neq 0$. Contradiction with Lemma
3.5.8 from \cite{unique}.
\end{proof}

The following notion plays a central role in the proof.
\begin{defi}\label{Def:2.3.3}
An element $D\in \D$  is said to be {\it hidden} if any
noninvertible function from $\K[X]^{(B)}$ is zero on $D$,
equivalently, $\langle \varphi_D,\mu\rangle>0$ for any $\mu\in
\X^+\setminus -\X^+$.
\end{defi}

The set of all hidden elements of $\D$ is denoted by
$\underline{\D}^\varnothing(=\underline{\D}^\varnothing_{G,X})$. Set
$\overline{\D}^\varnothing:=\D\setminus\underline{\D}^\varnothing$.
Let us establish some  properties of $\underline{\D}^\varnothing$
and $ \overline{\D}^\varnothing$.

\begin{Lem}\label{Lem:2.3.7}
 $\D^G\subset \overline{\D}^\varnothing$.
\end{Lem}
\begin{proof}
Let $D\in \D^G\cap \underline{\D}^\varnothing$. Since $G$ is not a
torus and $(G,G)$ acts nontrivially on $X$, it follows from
\cite{unique}, Lemma 3.5.7, that $\D\neq \{D\}$. Let $\lambda\in \X$
be such that $\langle \lambda,\varphi_{D'}\rangle\geqslant 0$ for
any $D'\in \D\setminus \{D\}$. Let us check that
$\ord_D(f_\lambda)\geqslant 0$. Assume the contrary. Choose a
function $g\in \K[X]^{(B)}$ that is zero on some $D'\in \D\setminus
\{D\}$. One can find positive $m,n$ such that $f^ng^m\in
\K[X]^{(B)}$ and $\ord_{D}(f^ng^m)=0$. Thence $f^ng^m$ is
invertible, which contradicts $\ord_{D'}(f)\geqslant
0,\ord_{D'}(g)>0$.  Therefore $\K[X]^{(B)}=\K[X\setminus D]^{(B)}$.
From  highest weight theory it follows that $\K[X]=\K[X\setminus
D]$, which is nonsense.
\end{proof}

\begin{Prop}\label{Prop:3.1.4}
Let  $\alpha,\alpha_1,\alpha_2$ be as in Lemma \ref{Lem:2.3.2}.
If $\D(\alpha_1)\subset\underline{\D}^\varnothing$, then
$\D=\D(\alpha_1)$ and $\varphi_{D_1}=\varphi_{D_2}$ for any
$D_1,D_2\in \D(\alpha_1)$ (the last condition is essential only if
$\alpha=\alpha_1$).
\end{Prop}
\begin{proof}
Consider the  case $\alpha=\alpha_1$. Let $D_1,D_2$ be different
elements of $\D(\alpha_1)$. Note that $\ord_{D_i}(f_\alpha)=1,
i=1,2$. By Lemma \ref{Lem:2.3.2}, $\ord_D(f_\alpha)\leqslant 0$ for
$D\in\D\setminus \D(\alpha_1)$. Choose $\lambda\in \X^+\setminus
-\X^+$. Set $m:=\min_{i=1,2}\langle\varphi_{D_i},\lambda\rangle$.
Since $\D(\alpha_1)\subset\underline{\D}^\varnothing$, we have
$m>0$. Put $f:=\frac{f_\lambda}{f_\alpha^m}$. Clearly,
$\ord_D(f)\geqslant 0$ for any $D\in \D$  and $\ord_{D_i}(f)=0$ for
some $i$. It follows that $f$ is an invertible element of
$\K[X]^{(B)}$. Therefore
$\langle\varphi_{D},\alpha\rangle=\langle\varphi_D,\lambda\rangle=0$
for $D\not\in \D(\alpha_1)$ and
$\langle\varphi_{D_1},\lambda\rangle=\langle\varphi_{D_2},\lambda\rangle$.
This implies the claim.

The proof in the remaining two cases is analogous (except for
$\alpha=\alpha_1+\alpha_2$ or $2\alpha_1$ one has to consider
$f=\frac{f_\lambda^2}{f_\alpha^m}$ instead of
$\frac{f_\lambda}{f_\alpha^m}$).
\end{proof}

\begin{Lem}\label{Lem:3.1.5}
In the notation of Proposition \ref{Prop:3.1.4}, $\Psi=\{\alpha\},
\Pi(\g)=\Supp(\alpha)$.
\end{Lem}
\begin{proof}
Consider the case $\alpha=k(\alpha_1+\alpha_2), k=1$ or
$\frac{1}{2}$ (the other two cases are analogous, even easier). Let
$\beta\in \Pi(\g)\setminus\Supp(\alpha)$ be such that
$\D(\beta)\neq\varnothing$. Corollary \ref{Cor:2.3.8} implies
$\D(\beta)\neq \D(\alpha_1)$, contradiction with Proposition
\ref{Prop:3.1.4}. So $\Pi(\g)\setminus \Supp(\alpha)\subset\Pi^a$.
By Lemma \ref{Lem:5.2}, $\Pi(\g)=\Supp(\alpha)$. So $[\g,\g]\cong \sl_2\times\sl_2$. Since $(\alpha_i,\alpha_1+\alpha_2)>0$ for $i=1,2$,
we see that $\Psi$ does not contain multiples of $\alpha_1,\alpha_2$.
Therefore $\Psi=\{\alpha\}$.
\end{proof}

\begin{Prop}\label{Prop:2.4.1}
Suppose $X\cong^G G*_HV$,
$\underline{\D}^\varnothing,\overline{\D}^\varnothing\neq\varnothing$.
Then $\underline{\D}^\varnothing$ consists of one element, say $D$,
and there exist a simple ideal $\g_1\subset\g, \g_1\cong\sl_n,$ and
$i\in \{1,n-1\}$ such that
\begin{enumerate}
\item The simple root $\alpha_i$ of $\g_1$ is the unique simple root of $\g$
positive on $\X^+\setminus -\X^+$.
\item $G_D=P_{\Pi(\g)\setminus\{\alpha_i\}}$.
\item $H$ is $G$-conjugate to a subgroup in $G_D$.
\end{enumerate}
\end{Prop}
\begin{proof}
We may assume that $G=Z(G)^\circ\times G_1\times\ldots\times G_k$,
where $G_i$ is a simple simply connected group. Let $H_i,
i=\overline{1,k},$ denote the projection of $H$ to $G_i$ and
$\rho_i$ denote the natural projection $G*_HV\twoheadrightarrow
G_i/H_i$.

{\it Step 1.} Since both
$\overline{\D}^\varnothing,\underline{\D}^\varnothing$ are nonempty, we see that no multiple of
any element in $\overline{\D}^\varnothing$ is principal (otherwise, there is a
noninvertible $B$-semiinvariant function nonvanishing on a hidden divisor). So $\Pic(X)$
is infinite.   Note that
$\Pic(X)\cong\Pic(G/H)\cong \X(H)/p(\X(G))$ (the last isomorphism is
due to Popov, \cite{Popov_Pic}), where $p$ is the restriction of
characters. It follows that $\h^H\cap[\g,\g]\neq \{0\}$. In
particular, for some $i\in \{1,\ldots,k\}$ the group $H_i$ is not semisimple
whence $H_i\neq G_i$. To
be definite, suppose $H_1\neq G_1$.

{\it Step 2.} Let $Y$ be an affine spherical $G$-variety of positive
dimension, $\varphi:X\rightarrow Y$ a  smooth $G$-morphism and $D\in
\underline{\D}^\varnothing$. If $Y\not\cong G/G_0$ with
$(G,G)\subset G_0$, then $D\in
\varphi^*(\underline{\D}_{G,Y}^\varnothing)$. Indeed, if
$\overline{\varphi(D)}=Y$, then $\ord_D(\varphi^*(f))=0$ for any $
f\in \K[Y]^{(B)}$. By the assumption on $Y$, there is a
noninvertible element in $\K[Y]^{(B)}$. Thus $\overline{\varphi(D)}$
is a divisor. To see that $D\in
\varphi^*(\underline{\D}_{G,Y}^\varnothing)$ note that
$\ord_{\overline{\varphi(D)}}(f)=\ord_D(\varphi^*(f))$. We also
remark that $G_{D}=G_{\overline{\varphi(D)}}$.

Applying this observation to $\rho_i:X\twoheadrightarrow G_i/H_i$ we
get that either
$\underline{\D}^\varnothing=\rho_i^*(\underline{\D}^\varnothing_{G,G_i/H_i})$
or $H_i=G_i$. It follows that $G_i=H_i$ for any $i>1$ and
$\underline{\D}^\varnothing_{G_1,G_1/H_1}\neq \varnothing$.

{\it Step 3.} Let us check that $\g_1\cong\sl_n$ and $H_1$ is
conjugate to a subgroup in $Z_{G_1}(\pi_1)$, where $\pi_1$ is the
fundamental weight of $\sl_n$ corresponding to the simple root
$\alpha_1$.

By step 1, $\h_1^{H_1}\neq \{0\}$. Put $L_1=Z_{G_1}(\h_1^{H_1})$.
Let us check now that $\rank_{G_1}(G_1/L_1)=1$. Assume the contrary.
Let $G_1/L_1$ be symmetric. By results of Vust,
\cite{Vust2}, $\varphi_D,D\in \D_{G,G_1/L_1},$ is the half of a
simple coroot of the symmetric space $G_1/L_1$. In particular, for
any $D\in \D_{G,G_1/L_1}$ there exists a noninvertible
(=nonconstant) function $f\in \K[G_1/L_1]^{(B)}$ such that
$\ord_D(f)=0$, so $D\in\overline{\D}_{G,G_1/L_1}^\varnothing$. 

So $G_1/L_1$ is not symmetric. We deduce from the classification
of \cite{Kramer} that the pair $(\g_1,\l_1)$ is either
$(\so_{2n+1},\gl_n)$ or $ (\sp_{2n},\sp_{2n-2}\times\K)$. However,
in both these cases $\l_1^{N_{G_1}(\l_1)}=0$. Applying step 2 to the
projection $G_1/H_1\rightarrow G_1/N_{G_1}(\l_1)$, we get
$\underline{\D}^\varnothing_{G,G_1/H_1}=\varnothing$. This
contradicts step 2.

So $\rank_{G_1}(G_1/L_1)=1$.  It follows from the classification in
\cite{Kramer} that $\g_1=\sl_n, \l_1\sim_G\z_{\g_1}(\pi_1)$. We note
that $\X^+_{G,G_1/L_1}=\Z_{\geqslant 0}(\pi_1+\pi_{n-1}),
\#\D_{G,G_1/L_1}=2$, the stabilizers of elements of
$\#\D_{G,G_1/L_1}$ are $P_{\Pi(\g)\setminus\{\alpha_i\}}, i=1,n-1$,
and the classes of two elements of $\rho_1^*(\D_{G,G_1/L_1})$ in
$\Pic(X)$ are opposite. By step 2,
$\#\underline{\D}^\varnothing\leqslant 2$.

Note that, by the above, any Levi subgroup of $G_1$ containing $H_1$ is
conjugate to $L_1$. Therefore $\dim \h_1^{H_1}=1$. Since $H_i=G_i$
for any $i>1$, it follows that $\rank\Pic(X)=1$.

{\it Step 4.} Let us check that $\#\underline{\D}^\varnothing=1$.
For $D\in \D$ we denote by $[D]$ the class of $D$ in
$F:=\Pic(X)/\operatorname{Tor}(\Pic(X))$. The group $F$ is generated
by $[D]$ for any $D\in \D_{G,G_1/L_1}$. The equality
$\#\underline{\D}^\varnothing=1$ will follow if we check that
 for any $D\in \underline{\D}^\varnothing, D'\in \D, D\neq D'$
 there is a positive integer $n$ such that $[D']=-n[D]$.  Assume the contrary. Then there is $f\in
\K(X)^{(B)}$ such that $(f)=D'-nD, n\geqslant 0$. If $n=0$, then $f$
is a noninvertible element of  $ \K[X]^{(B)}$ and $f$ is nonzero on
$D$. So let $n>0$. There is a function $g\in \K[X]^{(B)}$ such that
$\ord_D(g)=kn, k\in\N$. It remains to note that $f^{k}g\in
\K[X]^{(B)}, f^kg$ is noninvertible and $\ord_D(f^{k}g)=0$,
contradiction.

Without loss of generality, assume that
$\underline{\D}^\varnothing\subset \D(\alpha_1)$.

{\it Step 5.} If $\mu\in \X^+$ is such that 
$\langle\alpha_1^\vee,\mu\rangle=0$, then
$\langle\varphi_D,\mu\rangle=0$. Indeed, since
$\langle\alpha_1^\vee,\mu\rangle=0$, we see that $f_\mu\in
\K[X]^{(P_{\alpha_1})}$. To prove the claim  note that
$\overline{P_{\alpha_1} D}=X$ whence $f_\mu$ is nonzero on $D$.

{\it Step 6.} From step 5 it follows that
$\langle\alpha_1^\vee,\X^+\setminus -\X^+\rangle>0$. To prove the
proposition it remains to  check that there is no other simple root
with this property. Assume the contrary, let $\alpha\in
\Pi(\g)\setminus\{\alpha_1\}$ be such that $\langle\alpha^\vee,
\X^+\setminus -\X^+\rangle>0$. Since $\pi_1+\pi_{n-1}\in \X^+$, we
see that $\alpha=\alpha_{n-1}$ whence $n\geqslant 2$. Let $D_1$
denote the unique element in $\rho_1^*(\D_{G,G_1/L_1}(\alpha_{n-1}))$.

Suppose that there is  $D_2\in\D\setminus \D(\alpha_{n-1})\setminus
\{D\}$. Then, by step 4, there exists $\mu\in \X^+$ such that
$(f_\mu)=nD+D_2, n>0$. In particular, $f_\mu\in
\K[X]^{(P_{\alpha_{n-1}})}$ whence
$\langle\alpha_{n-1}^\vee,\mu\rangle=0$.

 It remains to consider the case when $\D=\{D\}\cup\D(\alpha_{n-1})$.
Suppose $\D(\alpha_{n-1})=\{D_1\}$. Recall that $F\cong \Z,
[D]+[D_1]=0$. It follows that $(f)$ is proportional to $D+D_1$ for
any $f\in \K(X)^{(B)}$. Therefore $D_1\in
\underline{\D}^\varnothing$, contradiction.

Suppose $\#\D(\alpha_{n-1})=2$. Let $D_2$ be the unique element in
$\D(\alpha_{n-1})\setminus\{D_1\}$. By Proposition \ref{Prop:2.3.1},
$\alpha_{n-1}\in \Psi$ and
$\langle\varphi_{D_i},\alpha_{n-1}\rangle=1,i=1,2$. By Lemma
\ref{Lem:2.3.2}, $\langle\varphi_D,\alpha_{n-1}\rangle\leqslant 0$.
It follows that $[D_1]+[D_2]=m[D]$ for some nonnegative $m$.
Contradiction with step 4.
\end{proof}

\begin{proof}[Proof of Theorem \ref{Thm:5.1}]
Below in the proof we write $\X^+,\X,\a,\Pi^a$ instead of
$\X^+_{G,X_i},\X_{G,X_i},\a_{G,X_i},$ $\Pi(\g)^a_{X_i}$ (recall that
$\Pi(\g)^a_{X_1}=\Pi(\g)^a_{X_2}$ in virtue of Corollary
\ref{Cor:2.2.4}), $\D_i,\underline{\D}_i^\varnothing,
\overline{\D}^\varnothing_i,\Psi_i,\Pi^b_i,\ldots,\Pi^d_i$ instead
of $\D_{G,X_i}$, etc.

Assume the contrary: there is no bijection $\iota:\D_1\rightarrow
\D_2$ with the desired properties. Thanks to the assumptions made in
the end of Section \ref{SECTION_monoid} and Corollary
\ref{Cor:2.2.3}, the theorem holds for $X_1(\mu),X_2(\mu)$ for any
$\mu\in \X^+\setminus -\X^+$, see Corollary \ref{Cor:2.2.3}.

It follows from Lemma \ref{Lem:1.1.1} that both actions
$G:X_1,G:X_2$ are locally effective. We may also assume that
$Z(G)^\circ$ acts on $X_1,X_2$ effectively. Thanks to Lemma
\ref{Lem:1.1.7} both $X_1,X_2$ are indecomposable.

{\it Step 1.} Here we construct a bijection
$\iota:\overline{\D}^\varnothing_1\rightarrow
\overline{\D}^\varnothing_2$ such that
$\varphi_{D}=\varphi_{\iota(D)}$. By our assumptions, for any
$\mu\in \X^+\setminus -\X^+$ there is a bijection
$\iota_{\mu}:\D_{G_\mu,X_1(\mu)}\rightarrow \D_{G_\mu,X_2(\mu)}$
such that $\varphi_D=\varphi_{\iota_\mu(D)}, G_D\cap
G_\mu=G_{\iota_\mu(D)}\cap G_\mu$. By Proposition~\ref{Prop:2.2.2},
there is a natural embedding $\D_{G_\mu,X_i(\mu)}\hookrightarrow
\D_{i}$ with image $\{D\in \D_i| \langle\mu,\varphi_D\rangle=0\}$.
In the sequel we identify $\D_{G_\mu,X_i(\mu)}$ with this image.
Since $\D_{G_\mu,X_i(\mu)}=\{D\in \D_i| f_\mu|_D\neq 0\}$,  we have
\begin{equation}\label{eq:2.3:1}\D_{G_\lambda,X_i(\lambda)}\cap
\D_{G_\mu,X_i(\mu)}=\D_{G_{\lambda+\mu},X_i(\lambda+\mu)},\forall
\lambda,\mu\in \X^+\setminus -\X^+.\end{equation} By definition,
\begin{equation}\label{eq:2.3:2}\overline{\D}^\varnothing_i=\bigcup_{\mu\in\X^+\setminus -\X^+}\D_{G_\mu,X_i(\mu)}.\end{equation}

We remark that $\iota_\lambda,\iota_\mu$, in general, do not
coincide on $\D_{G_{\lambda+\mu},X_1(\lambda+\mu)}$. Indeed, any
$G$-\! equivariant automorphism $\psi$ of $X_1$ induces a
bijection $\D_1\rightarrow \D_1, D\mapsto \psi(D)$. This bijection
can be nontrivial (take $\SL_2$ for $G$ and $G/T$ for $X_1$).

 Let $\mu_1,\ldots,\mu_k\in \X^+$ be the minimal set of generators of $\X^+$
modulo $X^+\cap -\X^+$, i.e., $\{\mu_1,\ldots,\mu_k\}=X^+\setminus \X^++\X^+$. 
For a subset $I\subset \{1,2,\ldots,k\}$ put $\mu_I:=\sum_{j\in
I}\mu_j,M_I:=G_{\mu_I}, \D^I_i:=\D_{G_{\mu_I},X_i(\mu_I)},$ $
\overline{\D}^I_i:=\cup_{J\supsetneq I}\D^J_i,
\underline{\D}^I_i:=\D^I_i\setminus\overline{\D}^I_i, i=1,2$. For
$I\neq\varnothing$ set $\iota_I:=\iota_{\mu_I}$. It follows from
(\ref{eq:2.3:1}) that $\D^{I}_j=\cap_{i\in I}\D_j^{\{i\}}$ whence
$\D_j^{I}\cap\D_j^{J}=\D_j^{I\cup J},j=1,2$. Note that
$\iota_I(\D_1^J)=\D_2^J$ for any $J\supset I$ because $\D_j^J=\{D\in
\D_j^I| \langle\varphi_D,\sum_{i\in J}\mu_i\rangle=0\}$ and
$\varphi_{\iota_I(D)}=\varphi_D$.  Thus the map
$$\iota:\overline{\D}^\varnothing_1\rightarrow \overline{\D}^\varnothing_{2},\iota|_{\underline{\D}^I_1}=\iota_I,$$ is a well-defined bijection
such that $\varphi_{\iota(D)}=\varphi_D$. Clearly, $\iota^{-1}$
coincides with $\iota^{-1}_I$ on $\underline{\D}^I_2$.

{\it Step 2.} Choose $\alpha_1\in \Pi(\g)\setminus \Pi^a$. Suppose
$\D_i(\alpha_1)\subset \underline{\D}_i^\varnothing$. Let us check
that $\alpha_1\in \Pi_{i}^d$ and
$G_D=P_{\Pi(\g)\setminus\{\alpha_1\}}$. To be definite, put $i=1$.
Assume the contrary. Then one of the following cases takes place.

{\it Case 1.} $\alpha_1\in \Pi_{1}^b$. Let $D_1^+,D_1^-$ denote
different elements of $\D_{1}(\alpha_1)$. By Proposition
\ref{Prop:3.1.4}, $\D_1=\D_{1}(\alpha_1)$ and
$\varphi_{D_1^+}=\varphi_{D_1^-}$.  By Lemma \ref{Lem:3.1.5},
$[\g,\g]=\sl_2$.  If $\Psi_{1}=\Psi_{2}$, then $\alpha_1\in
\Psi_{2}$. Applying Proposition \ref{Prop:3.1.4} again, we see that
$\varphi_{D_2^+}=\varphi_{D_2^-}$ for $D_2^+,D_2^-\in \D_{2}$. Set
$\iota(D_1^\pm)=D_2^\pm$.

Now suppose  $X_1,X_2$ are smooth.  Let $H\subset G,V$ be such that
$X_1\cong^G G*_HV$.  Let $\pi$ denote the natural projection $X_1\twoheadrightarrow G/H$. 
  By Proposition 4.2.4 from \cite{unique},
 the pull-back map
$\pi^*:\D_{G,G/H}\rightarrow\D_1$ is a bijection. Take $f_\lambda\in \K[X_1]^{(B)}_\lambda$.
The divisor of $f_\lambda$ is a pull-back of a divisor on $G/H$ whence $f_\lambda$ is constant 
on fibers of the vector bundle $\pi:X_1\rightarrow G/H$. This means that $f_\lambda\in \pi^*\K[G/H]$.
Therefore $V=\{0\}$ and $X_1=G/H$. 
 By Theorem 2
from \cite{unique}, $N_G(H)$ is not connected, whence $H=T_1\times
T_0$, where $T_1$ is a maximal torus of $\SL_2$, and $T_0\subset
Z(G)^\circ$. By our assumptions $X_1$ is indecomposable. Thus
$\g=\sl_2$. Now it is easy to check that $X_1\cong^G X_2$.

{\it Case 2.} $\alpha_1\in \Pi_{1}^c$. By Proposition
\ref{Prop:3.1.4}, it is enough to check that $\alpha_1\in \Pi_2^c$.
Similarly to the previous case, one can check that $G=\SL_2,
X_1=G/N_G(T)$. Then it is easy to see that $X_1\cong^G X_2$.

{\it Case 3.} $\alpha_1\in \Pi_1^d$ and there is $\alpha_2\in
\Pi(\g)\setminus \{\alpha_1\}$ such that
$\D_{1}(\alpha_1)=\D_{1}(\alpha_2)$. Again we only need to consider
the case when $X_1,X_2$ are smooth. Analogously to case 1,
$[\g,\g]\cong \sl_2\times\sl_2, X_1=G/H$. From
$\D_{1}(\alpha_1)=\D_1(\alpha_2)$ one can deduce  that $\h$ contains
the diagonally embedded subalgebra $\sl_2\subset [\g,\g]$.  Again,
since $X_1$ is indecomposable, we see that $\g=\sl_2\times\sl_2$.
Then we easily check that $X_1\cong^G X_2$.

{\it Step 3.} Let $D\in
\overline{\D}_{1}^\varnothing\cap\D_{1}(\alpha)$. In this step we
show that
\begin{enumerate}
\item If $\alpha\in \Pi_{1}^d$, then $\iota(D)\in
\D_{2}(\alpha)$ and $\alpha\in \Pi_{2}^d$.
\item If $\alpha\in \Pi_1^c$, then $\iota(D)\in
\D_{2}(\alpha)$ and $\alpha\in \Pi_{2}^c$.
\item If $\alpha\in \Pi_{1}^b$, then either $\iota(D)\in
\D_{2}(\alpha)$ and $\alpha\in \Pi_{2}^b$ or $\alpha\in \Pi_{2}^d$
and $\D_{2}(\alpha)\cap\overline{\D}^\varnothing_{2}=\varnothing$.
\end{enumerate}
 Let $I\subset \{1,\ldots,k\}$ be such that $D\in
\underline{\D}^I_1$.

{\it Case 1.} Suppose  $\alpha\in \Pi^d_{1}$. Since  $G_D\cap
M_I=G_{\iota(D)}\cap M_I$, the inclusion $\iota(D)\in
\D_{2}(\alpha)$ will follow if we check that $\alpha\in
\Delta(\m_I)$. Assume the contrary. Then  $f_{\mu_I}$ is not
${P_\alpha}$-semiinvariant. It follows that $(f_{\mu_I})$ contains a
$P_\alpha$-unstable prime divisor. But $D$ is the only such divisor.
Since $D\in \D_I$, we have $\ord_{D}(f_{\mu_I})=0$. Contradiction.
The inclusion $\alpha\in \Pi^d_{2}$ is easily deduced from
Proposition \ref{Prop:2.3.1} and the equality
$\varphi_{\iota(D)}=\varphi_D=\alpha^\vee|_\a$.

{\it Case 2.} Suppose $\alpha\in \Pi^c_{1}$. Then
$\varphi_{\iota(D)}=\varphi_D=\frac{\alpha^\vee}{2}|_\a$. If
$\alpha\in \Pi^d_{2}$, then $\D_{2}(\alpha)\subset
\overline{\D}_2^\varnothing$ because
$\varphi_{D'}=2\varphi_{\iota(D)}$ for $D'\in \D_2(\alpha)$. This
contradicts the previous case (recall that $X_1,X_2$ have equal
rights). Therefore $\alpha\in \Pi^b_{2}\cup\Pi^c_{2}$. By Corollary
\ref{Cor:2.3.8}, $\iota(D)\in \D_{2}(\alpha)$. If $\alpha\in
\Pi_{2}^b$, then there is $D'\in
\D_{2}(\alpha)\setminus\{\iota(D)\}$ such that
$\varphi_{D'}=\varphi_D$. Therefore $D'\in \underline{\D}_2^I$ and
$\iota^{-1}(D')\in \underline{\D}_1^I$. Again, by Corollary
\ref{Cor:2.3.8}, $\iota^{-1}(D')\in \D_{1}(\alpha)$. Contradiction
with $\alpha\in \Pi_1^c$.

{\it Case 3.} Suppose $\alpha\in \Pi_{1}^b$. If $\D_{2}(\alpha)\cap
\overline{\D}_2^\varnothing\neq \varnothing$, then, thanks to cases
1,2, $\alpha\in \Pi_{2}^b$. If
$\D_{2}(\alpha)\subset\underline{\D}_2^\varnothing$, then, by step
2, $\alpha\in \Pi_{2}^d$, q.e.d.

From case 2 it follows that $\Pi_1^c=\Pi_2^c$.

{\it Step 4.} Suppose $\Pi^b_{1}=\Pi^b_{2},
\Pi^d_{1}=\Pi^d_{2}$. 
Let us construct a bijection $\iota:\D_{1}\rightarrow \D_{2}$
extending the bijection $\overline{\D}_1^\varnothing\rightarrow
\overline{\D}_2^\varnothing$ constructed above and such that
$G_{\iota(D)}=G_D, \varphi_{\iota(D)}=\varphi_D$ for any $D\in
\D_1$.
 By step 3, for any $\alpha\in \Pi(\g), D\in \overline{\D}^\varnothing_{1}\cap\D_1(\alpha)$
 we have  $\iota(D)\in \D_2(\alpha)$.
 It follows that
$G_D=G_{\iota(D)}$.

Now choose $D\in \underline{\D}^\varnothing_1$.  Thanks to steps
2,3, exactly one of the following possibilities holds:
\begin{enumerate}
\item $\alpha\in \Pi_{1}^b$ and
$\overline{\D}_{i}^\varnothing\cap\D_{i}(\alpha)\neq \varnothing$
for any $\alpha\in \Pi(\g)$ such that $D\in \D_{1}(\alpha)$ and any
$i=1,2$.
\item There  is a unique root $\alpha\in \Pi(\g)$ such that $D\in
\D_{1}(\alpha)$.
This root is in $\Pi^d_1=\Pi^d_2$. Let $D'$ be the unique element of
$\D_2(\alpha)$. Then $G_D=G_{D'}=P_{\Pi(\g)\setminus\{\alpha\}},
D'\in \underline{\D}_2^\varnothing$.
\end{enumerate}
 In
case 2 put $\iota(D):=D'$. In case 1 set $A:=\{\alpha\in
\Pi(\g)|D\in \D_1(\alpha)\}$.  For $\alpha\in A$ let $D_\alpha$
denote the unique element in $\D_{1}(\alpha)\cap
\overline{\D}_1^\varnothing$. Clearly,
$\varphi_D\neq\varphi_{D_\alpha}$ for any $\alpha\in A$. Further,
$\langle\alpha,\varphi_{D_\beta}\rangle=\langle\alpha,\beta^\vee-\varphi_D\rangle<0$
whenever $\alpha,\beta$ are different elements of $A$. By step 3,
$\iota(D_\alpha)\in \D_{2}(\alpha)$. For $\alpha\in A$ let
$D_\alpha'$ denote the unique element in
$\D_{2}(\alpha)\setminus\{\iota(D_\alpha)\}$. Note that
$\langle\varphi_{D_\alpha'},\beta\rangle=\langle\alpha^\vee-\varphi_{D_\alpha},\beta\rangle=1$
for any $\alpha,\beta\in A$. Applying Proposition \ref{Prop:2.3.2}
and Corollary \ref{Cor:2.3.8}, we see that $D_\alpha'=D_\beta'$  for
any $\alpha,\beta\in A$. Put $\iota(D):=D'_\alpha,\alpha\in A$.

Since $\langle\varphi_D,\gamma\rangle\leqslant 0$ for any $\gamma\in
\Pi_{1}^b\setminus A$, we see that $\iota:\D_{1}\rightarrow\D_{2}$
is injective. As  already checked,  $G_D\subset G_{\iota(D)}$
for any $D\in \underline{\D}_1^\varnothing$. Finally, by the
symmetry between $X_1,X_2$, we obtain that $G_D=G_{\iota(D)}$ and
$\iota:\D_{1}\rightarrow \D_{2}$ is surjective.

{\it Step 5.} We may assume that $X_1,X_2$ are smooth and there is
$\alpha\in \Pi_{1}^d\cap\Pi_{2}^b$. In this case, according to case
3 of step 3, $\D_{1}(\alpha)\subset \underline{\D}_1^\varnothing$.

Since $\Pi_{2}^b\neq \varnothing$, it follows from step 2 that
$\overline{\D}^\varnothing_{1},\overline{\D}^\varnothing_{2}\neq\varnothing$.
Therefore $X_1$ satisfies  conditions (1)-(3) of  Proposition
\ref{Prop:2.4.1}. Let $\g_1$ be as in Proposition
\ref{Prop:2.4.1}. Note that $\alpha$ is the unique simple root of $\g$
such that $\langle\alpha^\vee,\X^+\setminus -\X^+\rangle>0$. Suppose
$\underline{\D}^\varnothing_2\neq\varnothing$. It follows that $X_2$
satisfies conditions (1)-(3) for the same simple ideal $\g_1$. So we
have $G_{D_1}=G_{D_2}$ for $D_i\in \underline{\D}^\varnothing_{i},
i=1,2$. It follows from Corollary \ref{Cor:2.2.3} and the
assumptions made in the end of Section \ref{SECTION_monoid} that
$X_1(D_1)\cong^M X_2(D_2)$, where $M$ is the standard Levi subgroup of
$G_{D_i}$. Suppose that $X_1\cong^G G*_{H_1}V_1, X_2\cong^G
G*_{H_2}V_2$. We may assume that $H_i\subset M, i=1,2$. Lemma
\ref{Lem:2.2.7} implies $X_i(D_i)\cong^{M} \Rad_u(\q^-)\times
M*_{H_i} V_i,i=1,2$, where $\q^-:=\b^-+\m$. By Remark
\ref{Rem:1.4.1}, $X_1\cong^G X_2$.

So it remains to consider the situation when
$\underline{\D}_2^\varnothing=\varnothing$. If  $\beta\in
\Pi^d_1\setminus\Pi^d_2, \beta\neq\alpha$, then, by step 3,
$\D_1(\beta)\in \underline{\D}^\varnothing_1$, which is impossible.
Analogously, $\Pi_2^d\subset \Pi_1^d$. So
$\Pi^b_{1}=\Pi^b_{2}\setminus \{\alpha\}$ and
$\Pi_{1}^d=\Pi^d_{2}\sqcup\{\alpha\}$. Note also that
$\D_{1}(\beta)\subset\overline{\D}_1^\varnothing$ for any $\beta\in
\Pi(\g),\beta\neq\alpha$.

Let $D_1$ denote the unique element of
 $\D_{1}(\alpha)$, $D^+_2,D^-_2$ denote  different elements of $\D_{2}(\alpha)$ and
$D_1^\pm:=\iota^{-1}(D_2^\pm)$. Then
$\varphi_{D_1}=\varphi_{D_1^+}+\varphi_{D_1^-}$. Put
$\D_1':=\D_{1}\setminus \{D_1,D^+_1,D_1^-\}, \D_2':=\D_{2}\setminus
\{D_2^+,D_2^-\}$. Let us check that $G_{\D'_1}=G_{\D_2'}$. Indeed,
choose $\beta\in \Pi(\g)$. Suppose there is $D'\in
\D_{2}(\beta)\cap\D_2'$. By step 3, $\iota^{-1}(D')\in
\D_{1}(\beta)$. Further, $\iota^{-1}(D')\neq D,D_1^\pm$.  Therefore
$G_{\D'_1}\subset G_{\D'_2}$. The opposite inclusion is proved in
the same way.

Let $M$ denote the standard Levi subgroup of $G_{\D'_i}$. Since
$\varphi_{D_1}=\varphi_{D_1^+}+\varphi_{D_1^-}$, we have
$$\X^+_{M,X_i(\D'_i)}=\{\lambda\in \X_{G,X}| \langle\lambda,\varphi_{D^\pm_j}\rangle\geqslant 0\}, j=1,2.$$
But  $\#\D_{M,X_1(\D_1')}\neq \#\D_{M,X_2(\D_2')}$. The assumptions made
in the end of Section \ref{SECTION_monoid} yield
$\D_i'=\varnothing$.

Since $\Pi_{2}^b\neq\varnothing$ and $\#\D_{2}=2$, we see that
$\Pi(\g)=\Pi^a\cup\Pi_{2}^b$. By Lemma \ref{Lem:5.2},
$[\g,\g]\cong\sl_2$. However, in this case Proposition
\ref{Prop:2.4.1} implies $\#\D_{1}(\alpha)=2$. Contradiction with
$\alpha\in \Pi_1^d$.
\end{proof}

Let us present an algorithm  recovering the set $\D_{G,X}$ with the
maps $D\mapsto \varphi_D, D\mapsto G_D$ from
$\X^+_{G,X},\Psi_{G,X}$.

\begin{Alg}\label{Alg:1}
Put $\X^+:=\X_{G,X}^+, \X:=\X_{G,X},\a:=\a_{G,X}, \Psi:=\Psi_{G,X},
\Pi^a:=\Pi(\g)^a_X,\ldots,\Pi^d:=\Pi(\g)^d_X$. We only need to
determine the set $\D'(G,\X^+,\Psi):=\D_{G,X}^G\cup_{\alpha\in
\Pi^b}\D_{G,X}(\alpha)$ and the maps $D\mapsto \varphi_D$ for $D\in
\D'(G,\X^+,\Psi)$. Then for $D\in \D'$ we recover $G_D$ from
$\varphi_D$ by using Corollary \ref{Cor:2.3.8}.

{\bf Step 1.} Compute $\{\mu_1,\ldots,\mu_k\}:=\X^+\setminus
(\X^+\cup\X^+)$. For $I\subset \{1,\ldots,k\}$ put
$M^I:=Z_G(\sum_{i\in I}\mu_i), \Pi^I=\{\alpha\in \Pi(\g)|
\langle\alpha^\vee, \sum_{i\in I}\mu_i\rangle=0\},
\X^{+I}:=\X^++\sum_{i\in I}\Z\mu_i, \Psi^I:=\{\alpha\in \Psi|
\Supp(\alpha)\subset \Pi^I\}$. Put $\D^I:=\D'(M^I,\X^{+I},\Psi^I),
\overline{\D}^I:=\cup_{J\supsetneq I}\D^J,
\underline{\D}^I:=\D^I\setminus \overline{\D}^I$. Note that
$\D'(G,\X^+,\Psi)=\sqcup_{I\subset \{1,2\ldots,k\}}\underline{\D}^I$
and $\underline{\D}^{\{1,\ldots,k\}}=\varnothing$.

{\bf Step 2$I$.} Here we compute the set $\underline{\D}^I$ together
with the map $D\mapsto \varphi_D$. Suppose we have already computed
the set $\overline{\D}^I$  together with the map $D\mapsto
\varphi_D$.  Note that
$\underline{\D}^I\subset\underline{\D}^\varnothing_{M_I,X(\sum_{i\in
I}\mu_i)}$.

{\it Case 1.} $\Pi^I=\Pi^a$.

{\it Case 1a.} $\X^{+I}=-\X^{+I}$. In this case
$\underline{\D}^I=\varnothing$.

{\it Case 1b.} $\rank(\X^{+I}\cap -\X^{+I})\leqslant \rank\X-2$.
Using the argument of the proof of Lemma \ref{Lem:2.3.7}, one can
show that $\underline{\D}^I=\varnothing$.

{\it Case 1c.} $\rank(\X^{+I}\cap -\X^{+I})=\rank\X-1$. Let $i$ be
such that the image of $\mu_i$ generates
$\X^{+I}/(\X^{+I}\cap-\X^{+I})$. Then $\underline{\D}^I=\{D\}$ and
$\varphi_D$ is given by $\langle\varphi_D,\mu_j\rangle=\delta_{ij},
\langle\varphi_D,\X^{+I}\cap -\X^{+I}\rangle=0$.

{\it Case 2.} $\Pi^I=\{\alpha\}\sqcup \Pi^a,\alpha\in \Pi^b,
\langle\alpha, \varphi_D\rangle\leqslant 0$ for any $D\in
\overline{\D}^I$. Then $\underline{\D}_I=\{D_1,D_2\}$ and
$\varphi_{D_i}=\frac{1}{2}\alpha^\vee|_{\a}$.

{\it Case 3.} Otherwise, let $\Sigma$ denote the set of all
$\alpha\in \Pi^I\cap \Pi^b\setminus \cup_{J\supsetneq I}\Pi^J$ such
that there is only one divisor $D\in \overline{\D}^I$ with
$\langle\alpha,\varphi_D\rangle=1$ and put $
\varphi_\alpha=\alpha^\vee-\varphi_{D}$. Let
$\varphi_1,\ldots,\varphi_l$ be all different values of
$\varphi_\alpha$. Then $\underline{\D}^I:=\{D_1,\ldots,D_l\}$ with
$\varphi_{D_i}=\varphi_i, i=\overline{1,l}$.
\end{Alg}

Note that in the algorithm one does not need to know  the whole set
$\Psi_{G,X}$ but only the subset consisting of  all roots $\alpha$
of one of the three forms indicated in Lemma \ref{Lem:2.3.2}.

\section{Equality of the systems of spherical
roots}\label{SECTION_roots} The goal of this section is to prove
Theorem \ref{Thm:6.1}. The proof will be given at the end of the
section. It is rather close in spirit to that of Theorem
\ref{Thm:5.1}. It is based on Proposition \ref{Prop:2.2.1} and
Theorem \ref{Thm:5.1} and uses the notion of a hidden spherical
root, see Definition \ref{defi:3.1.1}. We shall see that this proof is much
easier than that of Theorem \ref{Thm:5.1}.

Let $X$ be a spherical $G$-variety. We suppose that $G$ is not a
torus and the action $G:X$ is locally effective. $\a,\X,\X^+$, etc.
have the same meaning as in the previous section.

\begin{defi}\label{defi:3.1.1}
An element $\gamma\in \Psi$ is called {\it hidden} if the following
two conditions are satisfied.
\begin{enumerate}
\item For any  $D\in \D$ there is $\alpha\in
\Supp(\gamma)$ such that $D\in \D(\alpha)$.
\item $\gamma$ is not of any of the types (1)-(3) indicated in Lemma
\ref{Lem:2.3.2}.
\end{enumerate}
By $\Psi^0$ we denote the subset of $\Psi$ consisting of all
nonhidden roots.
\end{defi}

The  terminology is justified by step 1 of the proof of Theorem
\ref{Thm:6.1} below.

\begin{Prop}\label{Prop:3.1.3}
Suppose $\Psi\neq \Psi^0$. If $\#\Psi=1$, then $[\g,\g]$ is simple
and $\Supp(\gamma)=\Pi(\g)$ for $\gamma\in\Psi$. If $\#\Psi>1$, then
$(\Pi(\g),\Psi,\Pi^a)$ is one of the triples listed below:
\begin{enumerate}
\item $\Pi(\g)=C_n,n\geqslant 2,\Psi=\{k\alpha_1,\alpha_1+\alpha_n+2\sum_{i=2}^{n-1}\alpha_i\}$
for $k=1$ or $2$, $\Pi^a:=\{\alpha_3,\ldots,\alpha_n\}$.
\item $\Pi(\g)=G_2, \Psi=\{\alpha_2,\alpha_1+\alpha_2\},\Pi^a=\varnothing$.
\item $\Pi(\g)=C_n\times
A_1,n\geqslant
2,\Psi=\{\alpha_1+\alpha_1',\alpha_1+\alpha_n+2\sum_{i=2}^{n-1}\alpha_i\}$
(here $\alpha_1,\ldots,\alpha_n$,resp. $\alpha_1'$, are simple roots
in $C_n$, resp. $A_1$), $\Pi^a=\{\alpha_3,\ldots,\alpha_n\}$.
\item $\Pi(\g)=B_4, \Psi=\{\alpha_2+2\alpha_3+3\alpha_4,
\alpha_1+\alpha_2+\alpha_3+\alpha_4\},\Pi^a=\{\alpha_2,\alpha_3\}$.
\end{enumerate}
In all cases the second root is hidden.
\end{Prop}
\begin{proof}
Choose $\gamma\in \Psi\setminus \Psi^0$. Inspecting Table 1 from
\cite{Wasserman}, we see that $\Supp(\gamma)$ is connected. Recall
that $\langle\alpha^\vee,\a\rangle=0$ for any $\alpha\in \Pi^a$
(\cite{unique}, Lemma 3.5.8). In particular, if $\beta\in \Psi$ and
$\alpha\in \Pi^a$, then either $\alpha\in \Supp(\beta)$ or $\alpha$
is not adjacent to $\Supp(\beta)$.

{\it Step 1.} Let $\alpha\in \Pi^c$. We claim that $\alpha\in \Supp
\gamma$. Indeed, by Corollary \ref{Cor:2.3.8},
$G_D=P_{\Pi(\g)\setminus\{\alpha\}}$ for $D\in \D(\alpha)$. Further,
any simple root adjacent to $\alpha$ is not of type a) and
$\langle\alpha^\vee,\gamma\rangle\leqslant 0$.

{\it Step 2.} Let $\alpha\in \Pi^b$. Let $D_1,D_2$ be different
elements of $\D(\alpha)$. If $\beta$ is adjacent to $\alpha$, then
$\beta\not\in \Pi^a$. Let $\alpha\not\in\Supp(\gamma)$. Note that
$\langle\alpha^\vee,\gamma\rangle\leqslant 0$. Then, thanks to
Corollary \ref{Cor:2.3.8}, there are $\alpha^1,\alpha^2\in
\Pi^b\cap\Supp(\gamma)$ such that $\alpha^1\neq\alpha^2, D_i\in
\D(\alpha)\cap\D(\alpha^i),i=1,2$.

{\it Step 3.} Let $\alpha\in \Pi^d$. If
$\alpha\not\in\Supp(\gamma)$, then, by Corollary \ref{Cor:2.3.8},
there is (uniquely determined) $\alpha'\in \Pi^d\cap\Supp(\gamma)$
such that $\D(\alpha)=\D(\alpha'),
\langle\alpha^\vee-\alpha'^\vee,\a\rangle=0$ and
$k(\alpha+\alpha')\in \Psi$ for $k=1$ or $\frac{1}{2}$. Again,
$\beta\not\in \Pi^a$ whenever $\beta$ is adjacent to $\alpha$ or
$\alpha'$. Since $\langle\alpha^\vee,\gamma\rangle\leqslant 0$, we
have $\langle\alpha'^\vee,\gamma\rangle\leqslant 0$.

{\it Step 4.} If $\#\Psi=1$, then, by the previous steps,
$\Supp(\gamma)\cup\Pi^a=\Pi(\g)$. Applying Lemma \ref{Lem:5.2}, we
see that $\Supp(\gamma)=\Pi(\g)$. Till the end of the proof we
suppose that $\#\Psi>1$. Exactly one of following possibilities
takes place:
\begin{itemize}
\item[(A)] $\Pi^a\cup\Pi^d=\Pi(\g)=\Supp(\gamma)$.
\item[(B)] There are two adjacent roots $\alpha^1,\alpha^2\in
\Supp(\gamma)\setminus \Pi^a$ such that
$\langle\alpha^{1\vee},\gamma\rangle\leqslant 0$ and either
$\alpha^1\in \Pi^b\cup\Pi^c$ or $\alpha^1\in \Pi^d$ and there is
$\beta\in \Pi^d$ such that $\D(\beta)=\D(\alpha^1)$.
\end{itemize}

{\it Step 5.} Let us consider possibility (B). Inspecting Table 1 in
\cite{Wasserman}, we see that only the following three cases are
possible:
\begin{itemize}
\item[(BC)] $\Supp(\gamma)=C_n, n\geqslant 2,
\gamma=\alpha_1+\alpha_n+2\sum_{i=2}^{n-1}\alpha_i,
\alpha^1=\alpha_1,\alpha^2=\alpha_2$.
\item[(BG)] $\Supp(\gamma)=G_2, \gamma=\alpha_1+\alpha_2, \alpha^1=\alpha_2,\alpha^2=\alpha_1$.
\end{itemize}
If $\Pi^b\not\subset \Supp(\gamma)$, then, according to step 2,
$\#\Pi^b\cap\Supp(\gamma)\geqslant 2$. Note also that
$\langle\beta,\gamma\rangle\leqslant 0$ for any $\beta\in \Pi^b$.
Now we check case by case that $\Pi^b\subset \Supp(\gamma)$.

 Consider case (BG). Since $\langle\alpha^{1\vee},\gamma\rangle<0, \langle\beta^\vee,\gamma\rangle=0$
for any $\beta\in \Pi(\g)\setminus \Supp(\gamma)$, we get
$\alpha^1\in \Pi^b$. Thus $\Pi(\g)=\Supp(\gamma)\cup\Pi^a$. By Lemma
\ref{Lem:5.2}, $\Pi(\g)=G_2$. So we get possibility 2 of the
proposition.

Consider case (BC). Note that
$\Supp(\gamma)\setminus\{\alpha_1,\alpha_2\}\subset\Pi^a,
\langle\alpha_1^\vee,\gamma\rangle=0$. Suppose $\alpha_1\in
\Pi^b\cup\Pi^c$. Since $\langle\alpha_2,\gamma\rangle>0$, we get
from step 3 that  $\Pi^d\subset \Supp(\gamma)$.  So
$\Pi(\g)\setminus \Supp(\gamma)\subset \Pi^a$. By Lemma
\ref{Lem:5.2}, $\Pi(\g)=\Supp(\gamma)$ and we get possibility 1 of
the proposition.

Now suppose $\alpha_1\in \Pi^d$. Let $\beta$ be a (unique, see
Corollary \ref{Cor:2.3.8}) element  of $\Pi^d\setminus\Supp(\gamma)$
such that $\D(\beta)=\D(\alpha_1)$. Since $\langle\alpha_1^\vee,
\gamma\rangle=0$, we see that $\beta$ is not adjacent to
$\Supp(\gamma)$. For $k=1$ or $\frac{1}{2}$ we get
$\gamma_1:=k(\alpha_1+\beta)\in \Psi$. Since $\alpha_2\in \Pi^d$, we
have $k=1$.  By step 3, $\Pi^d\setminus\Supp(\gamma)=\{\beta\}$. It
follows that $\Supp(\gamma)\cup\Supp(\gamma_1)\cup\Pi^a=\Pi(\g)$
whence $\Pi(\g)=\Supp(\gamma)\cup\{\beta\}$. We get possibility 2 of
the proposition.

{\it Step 9.} It remains to consider case (A).  Choose $\gamma_1\in
\Psi\setminus\{\gamma\}$. Recall that $\Psi\subset \Pi^{a\perp}$.
Therefore $\#\Pi^d\geqslant 2$. Inspecting Table 1 from
\cite{Wasserman}, we see that one of the following possibilities
takes place:
\begin{itemize}
\item[(AA)] $\Pi(\g)=A_n,n\geqslant 2, \gamma=\alpha_1+\ldots+\alpha_n,
\Pi^d=\{\alpha_1,\alpha_n\}$.
\item[(AB)] $\Pi(\g)=B_n,n\geqslant 2, \gamma=\alpha_1+\ldots+\alpha_n,
\Pi^d=\{\alpha_1,\alpha_n\}$.
\item[(AC)]  $\Pi(\g)=C_n, n> 2,
\gamma=\alpha_1+\alpha_n+2\sum_{i=2}^{n-1}\alpha_i,
\Pi^d=\{\alpha_1,\alpha_2\}$.
\item[(AG)] $\Pi(\g)=G_2, \gamma=\alpha_1+\alpha_2$.
\end{itemize}

Inspecting Table 1 from \cite{Wasserman} again and taking into
account that $(\gamma,\gamma_1)\leqslant 0$ and
$\langle\gamma_1,\alpha^\vee\rangle=0$ for any $\alpha\in \Pi^a$, we
get possibility 4 of the proposition.
\end{proof}

\begin{proof}[Proof of Theorem~\ref{Thm:6.1}]
It follows from Theorem \ref{Thm:5.1} that the type of a root
$\alpha\in \Pi(\g)$ is the same for $X_1$ and $X_2$. We put
$\Pi^a:=\Pi(\g)^a_{X_i}, \ldots,\Pi^d:=\Pi(\g)_{X_i}^d,
\a:=\a_{G,X_i},\X^+:=\X^+_{G,X_i}, \Psi_i:=\Psi_{G,X_i},i=1,2$. We
also identify $\D_{G,X_1},\D_{G,X_2}$ and write $\D$ instead of
$\D_{G,X_i}$.

Suppose $\Psi_1\neq\Psi_2$. Again, we impose the assumptions
in the end of Section \ref{SECTION_monoid}. In
particular, it follows from Corollary \ref{Cor:2.2.3} and Theorem
\ref{Thm:5.1} that $\Psi_{M,X_1(D)}=\Psi_{M,X_2(D)}$ for any $D\in
\D$, where $M$ denotes the standard Levi subgroup of $G_D$. So we
may assume that $\D^G=\varnothing$. As in the proof of Theorem
\ref{Thm:5.1} we may also assume that both actions $G:X_1,G:X_2$ are
locally effective and both $G$-varieties $X_1,X_2$ are
indecomposable.

{\it Step 1.} Let us show that $\Psi_1^0=\Psi_2^0$. It is enough to
check that $\Psi_1^0\subset \Psi_2^0$.

Let $\alpha\in \Psi_1^0$. If $\alpha$ has one of the forms indicated
in Lemma \ref{Lem:2.3.2}, then the inclusion $\alpha\in \Psi_2$
follows from Theorem \ref{Thm:5.1}, Propositions
\ref{Prop:2.3.1},\ref{Prop:2.3.2}.

Now let $D\in \D$ and $M$ be the standard Levi subgroup of $G_D$.
Suppose $\Supp\alpha\subset \Pi(\m)$. By Proposition
\ref{Prop:2.2.1}, $\alpha\in \Psi_{M,X_1(D)}=\Psi_{M,X_2(D)}\subset
\Psi_2$.

{\it Step 2.} So we may assume that $\Psi_1\neq\Psi_1^0$ whence
$\Psi_1$ is one of the systems described in Proposition
\ref{Prop:3.1.3}. Let us check that $\Psi_2=\Psi_1^0$. Otherwise
$\Psi_2^0\neq\Psi_2$ and $\Psi_2$ is also one of the systems from
Proposition \ref{Prop:3.1.3}. If $\#\Psi_1^0\neq\varnothing$, then
we get $\Psi_1=\Psi_2$ because all systems $\Pi(\g)$ in the list of
Proposition \ref{Prop:3.1.3} are distinct. So $\#\Psi_{i}=1$ and
$\Supp(\gamma_i)=\Pi(\g)$ for the unique element $\gamma_i\in
\Psi_{i}$. Using Table 1 from \cite{Wasserman} and the equality
$\Pi(\g)_{X_1}^a=\Pi(\g)_{X_2}^a$, we get $\gamma_1=\gamma_2$.

{\it Step 3.} Let us check that $X_2\cong^G G*_HV$, where
$(G,G)\subset H$. Assume the contrary, let $(G,G)\not\subset H$ or,
equivalently, $\widetilde{H}:=N_G(H)^\circ\neq G$. Let $H_0$ denote
the stabilizer of a point from the open $H$-orbit in $V$. Then
$H_0\subset \widetilde{H}$. It is clear that
$\Psi_{G,G/H_0}=\Psi_2$. Applying Proposition 3.4.3 from
\cite{unique} to the pair $H_0\subset \widetilde{H}$, we see that
$\Psi_{G,G/\widetilde{H}}=\varnothing$ or
$\Psi_{G,G/\widetilde{H}}=\Psi_{G,G/H_0}$. However,
$\a_{G,G/\widetilde{H}}$ is generated by $\Psi_{G,G/\widetilde{H}}$
(it follows, for example, from \cite{unique}, Lemma 3.1.4) and
$\a_{G,G/\widetilde{H}}$ contains a dominant weight. On the
contrary, $\Span_\Q(\Psi^0_1)$ does not contain a dominant weight.

{\it Step 4.} Suppose $\#\Psi_1=1$. Then $\Psi_2=\varnothing$. It is
known, see, for example, \cite{Knop5}, Corollary 6.2, that $H_v$
contains a maximal unipotent subgroup of $G$ for any $v\in V$. One
easily deduces from this that $V$ as a $(G,G)$-module is the
tautological $\SL_n$- or $\Sp_{2n}$-module. It follows that $\a\cap
[\g,\g]\subset \Q\pi_1$. But, according to Table 1 from
\cite{Wasserman}, $\Psi_1\not\subset \Q\pi_1$, contradiction.

Suppose  $\#\Psi_1=2$. Note that $(\Pi(\h),\Pi(\h)_V^a)=(\Pi(\g),\Pi(\g)^a_V)$ 
and so the l.h.s. is one of the four pairs listed in Proposition \ref{Prop:3.1.3}. Since
$\D^G=\varnothing$, we get $\D_{H^\circ,V}^{H^\circ}=\varnothing$.
The classification of spherical modules, see, for example,
\cite{Leahy}, shows that there are no pairs $(H^\circ,V)$, where $V$ is a
spherical $H^\circ$-module with 
 $\Pi(\h),\Pi(\h)^a_V$ listed in Proposition \ref{Prop:3.1.3} and 
$\D_{H^\circ,V}^{H^\circ}=\varnothing$. The set $\Pi(\h)^a_V$ is determined from Leahy's
tables as follows: this is the set of all simple roots annihilated by all highest weights 
in $\K[V]$. 
\end{proof}

\section{Invariant K\"{a}hler structures}\label{SECTION_Kahler}
In this section $K$ is a compact connected Lie group, $G$ is the
complexification of $K$, and $M$ is a multiplicity free compact
Hamiltonian $K$-manifold (see Introduction) with symplectic form
$\omega$ and moment map $\mu$. Put $\widetilde{\omega}=\omega+\mu$.
 This is an equivariantly closed form on $M$ called the {\it equivariant} symplectic
 form. We say that $\widetilde{\omega}$ is an {\it equivariant K\"{a}hler} form if
 $\omega$ is K\"{a}hler. We denote by
$[\widetilde{\omega}]$ the class of $\widetilde{\omega}$ in the
second equivariant cohomology group $H^2_K(M,\R)$.

As above, we fix a Borel subgroup $B\subset G$ and a maximal torus
$T\subset B$. We may assume that $T_K:=T\cap K$ is a maximal torus in
$K$. The choice of $B$ and $T$ defines the Weyl chamber $\t_+\subset
\k^*$. Define the invariant moment map $\psi:M\rightarrow \t_+$  as
in Introduction.

\begin{defi}\label{defi:2} We say that a complex structure $I$ on
$M$ is compatible (with $K,\omega$) if $I$ is $K$-invariant, and
$\omega$ is a K\"{a}hler form with respect to $I$.
\end{defi}

\begin{Prop}[\cite{Woodward_Kahler}, Proposition 5.2]\label{Prop:8.1}
Let $I$ be a compatible complex structure on $M$. Then the $K$-action
on $M$  extends to a unique action $G:M$ by holomorphic
automorphisms.  Moreover, $M$ is a spherical (algebraic) projective
$G$-variety.
\end{Prop}

This proposition allows one to define the valuation cone of $(M,I)$,
which we denote by $\V(M,I)$. The objective of this section is to
prove the following uniqueness result.

\begin{Thm}\label{Thm:8.1}
Let $I_1,I_2$ be two compatible complex structures on $M$. Suppose
$\V(M,I_1)=\V(M,I_2)$. Then there is a $K$-equivariant
diffeomorphism $\varphi:M\rightarrow M$ preserving
$[\widetilde{\omega}]$ and such that $\varphi^*(I_2)=I_1$.
\end{Thm}

The restriction $\V(M,I_1)=\V(M,I_2)$ is essential, see
\cite{Woodward_Kahler}, Remark 4.4.

To prove the theorem we need to recall some more or less standard
facts.

Let $X$ be a smooth projective $G$-variety. Denote by $\Pic_G(X)$
the equivariant Picard group of $X$. Choose $\L\in \Pic_G(X)$.
Suppose that $\L$ is very ample (as a usual line bundle). To $\L$
one can assign an equivariant K\"{a}hler form
$\widetilde{\omega}_{\L}$ as follows. Let $V$ denote the $G$-module
$H^0(X,\L)^*$ and $\iota:X\hookrightarrow \P(V)$ be the embedding
induced by $\L$. Choose a $K$-invariant hermitian form
$(\cdot,\cdot)$ on $V$.  Let $\omega_{FS}$ denote the corresponding
Fubini-Study form on $\P(V)$ and $\mu_{FS}$ be the corresponding
moment map:
$$\langle\mu_{FS}(x),\xi\rangle=\frac{(\xi v,v)}{2\pi
i(v,v)},$$ where $v$ denotes a nonzero vector on the line $x$.  Put
$\omega_{\L}:=\iota^*(\omega_{FS}), \mu_{\L}:=\mu_{FS}\circ\iota,
\widetilde{\omega}_\L:=\omega_L+\mu_L$. By $\psi_\L$ we denote the
invariant moment map associated with $\mu_\L$. Note that
$\widetilde{\omega}_\L$ does not depend (up to a $K$-equivariant
diffeomorphism) from the choice of $(\cdot,\cdot)$

 We have a unique homomorphism
$\Pic_G(X)\otimes_\Z\R\rightarrow H^2_K(X,\R)$ mapping a very
ample $G$-bundle $\L$ to the class $[\widetilde{\omega}_\L]$.

\begin{Lem}\label{Lem:8.4}
Suppose $X$ is  spherical. Then the  homomorphism
$\Pic_G(X)\otimes_\Z\R\rightarrow H^2_K(X,\R)$ is an isomorphism.
\end{Lem}
\begin{proof}
For $\chi\in \X(G)$ let $\C_\chi$ denote the trivial bundle on $X$
on which $G$ acts by  $\chi$. Note that
$\widetilde{\omega}_{\L\otimes
\C_\chi}=\widetilde{\omega}_{\L}+i\chi$. So we have a commutative
diagram, where the top sequence is exact

\begin{picture}(110,30)
\put(2,22){$0$}\put(12,22){$\X(G)\otimes_\Z\R$}\put(42,22){$\Pic_G(X)\otimes_\Z\R$}
\put(77,22){$\Pic(X)\otimes_\Z\R$}\put(107,22){$0$}
\put(2,2){$0$}\put(15,2){$(\k/[\k,\k])^*$}\put(47,2){$H^2_K(X,\R)$}\put(82,2){$H^2(X,\R)$}\put(107,2){$0$}
\put(5,23){\vector(1,0){6}}\put(32,23){\vector(1,0){9}}
\put(68,23){\vector(1,0){8}}\put(100,23){\vector(1,0){6}}
\put(5,3){\vector(1,0){10}}\put(32,3){\vector(1,0){14}}
\put(65,3){\vector(1,0){15}}\put(100,3){\vector(1,0){6}}
\put(24,20){\vector(0,-1){13}} \put(55,20){\vector(0,-1){13}}
\put(88,20){\vector(0,-1){13}}
\end{picture}

 As was noted in the proof of Proposition 5.2 in
\cite{Woodward_Kahler}, there is an action of $\C^\times$ on $X$
with finitely many fixed points. The Bialynicki-Birula decomposition
of $X$ induced by this action (see \cite{BB}) consists of affine
spaces. It follows that $H^1(X,\R)=\{0\}$ and the right vertical
arrow is an isomorphism. Since $H^1(X,\R)=\{0\}$, we see that the
the bottom sequence is exact. Note that the left vertical arrow is
an isomorphism. Since the top sequence is exact, we see that the
middle vertical arrow is an isomorphism.
\end{proof}

\begin{Lem}\label{Lem:8.5}
Let $X$ be a smooth projective $G$-variety. Then the following
assertions hold.
\begin{enumerate}
\item The subset  $H^2_K(X,\R)^+\subset H^2_K(X,\R)$ consisting of all classes
of equivariant K\"{a}hler forms is open.
\item The moment polytope for an equivariant K\"{a}hler form $\widetilde{\omega}$
depends only on $[\widetilde{\omega}]$ and this dependence is
continuous.
\end{enumerate}
\end{Lem}
\begin{proof}
By definition, $\omega$ is K\"{a}hler iff $\omega_x(iv,v)>0$ for all
$x\in X, v\in T_xX$ whence assertion 1. To prove the first claim of
assertion 2 one uses Moser's trick exactly as in the proof of
Proposition 5.2 in \cite{Woodward_Kahler}. The second claim stems
from the formula $\langle\mu(x),
[\xi,\eta]\rangle=\omega_x(\xi_*,\eta_*)$.
\end{proof}

\begin{Prop}\label{Prop:8.6}
Let $X$ be a smooth projective spherical $G$-variety and $\L$ a very
ample $G$-bundle. We consider $X$ as a Hamiltonian $K$-manifold with
respect to the equivariant form $\widetilde{\omega}_\L$. Then the
following assertions hold:
\begin{enumerate}
\item $2\pi i\im\psi_\L$ is a rational polytope and $(2\pi i\im\psi_\L)\cap\t^*(\Q)=
\bigcup_{d\in \N}\{\frac{\lambda}{d}| H^0(X,\L^{\otimes
d})^{(B)}_\lambda\neq \{0\}\}$.
\item Let $\sigma$ be a rational $B$-semiinvariant section of $\L$ of weight $\mu$.
Then
$$2\pi i\im\psi_\L=\mu+\{\lambda\in \a_{G,X}(\R)|\langle\lambda,\varphi_D\rangle\geqslant -\ord_D(\sigma),
\forall D\in\D_{G,X}\}.$$
\item Let $K_0$ denote the principal isotropy subgroup for the Hamiltonian
action $K:X$ (see Introduction). The group $\X(T_K/(T_K\cap K_0))$
coincides with $\X_{G,X}$.
\end{enumerate}
\end{Prop}
\begin{proof}
The first assertion is due to Brion, \cite{Brion_polytope}.
Assertion 2 easily follows from the first one. The third assertion
seems to be known but we failed to find its proof in the literature.
So we give a proof here.

Below we put $\omega:=\omega_\L, \mu:=\mu_\L$.

{\it Step 1.} Put $V=H^0(X,\L)^*$ and let $(\cdot,\cdot)$ be the
hermitian form on $V$ used to define $\widetilde{\omega}_\L$. Put
$\widetilde{K}:=K\times Z$, where $Z$ is a one-dimensional compact
torus, and let $Z$ act on $V$ by scalar multiplications. Denote by
$\widetilde{X}$ the affine cone over $X$. The map $\Phi:V\rightarrow
\widetilde{\k}^*$ given by $\langle\Phi(v),
\xi\rangle=\frac{1}{2i}(\xi v,v)$ is a moment map for the action
$\widetilde{K}:V$. Note that $\mu^{-1}(\t_+)$ is just the set of
lines in $\Phi^{-1}(\t_+\times\z)\cap\widetilde{X}$. Let
$\widetilde{K}_0$ denote the stabilizer of a point $v\in
\Phi^{-1}(\t_+\times\z)\cap \widetilde{X}$ in general position.  It
is easy to see that the restriction of the projection
$\widetilde{K}\rightarrow K$ to $\widetilde{K}_0$  is injective and
its image coincides with $K_0$. So it is enough to prove that
\begin{equation}\label{eq:8.2}\X(\widetilde{T}_K/(\widetilde{T}_K\cap\widetilde{K}_0))=\X_{\widetilde{G},\widetilde{X}},\end{equation}
where $\widetilde{T}_K:=T_K\times Z,
\widetilde{G}:=G\times\C^\times$.

{\it Step 2.} Let $\t_+^0$ denote the interior of the smallest face
of $\t_+$ containing $i\im\psi_\L$. Put
$\widetilde{Y}:=\Phi^{-1}(\t_+^0\times\z)$. Replacing $K$ with a
covering we may assume that $K$ is the direct product of a torus and
a simply connected semisimple group. It follows now from results of
\cite{GSJ} (see especially Theorem 4.9, Theorem 6.11, and the proof
of Theorem 6.17) that there is a $\widetilde{T}_K$-equivariant
embedding $\widetilde{Y}\rightarrow \widetilde{X}\quo
U:=\Spec(\C[\widetilde{X}]^U)$ with  dense image. To check
(\ref{eq:8.2}) note that
$\X_{\widetilde{G},\widetilde{X}}=\X(\widetilde{T}/\widetilde{T}_0)$,
where $\widetilde{T}:=T\times\C^\times$ and $\widetilde{T}_0$ is the
kernel for the action $\widetilde{T}:\widetilde{X}\quo
U$ (i.e. the kernel of the corresponding homomorphism $\widetilde{T}\rightarrow \Aut(\widetilde{X}\quo
U)$).
\end{proof}

Let us generalize  assertion 2 of Proposition \ref{Prop:8.6} to the
case of an arbitrary equivariant K\"{a}hler form
$\widetilde{\omega}$.
 Let $\psi$ denote the  invariant moment map corresponding to $\widetilde{\omega}$.

Consider the embedding $\X_{G,X}\rightarrow \X(T)\times
\Z^{\D_{G,X}}$ given by $\lambda\hookrightarrow(\lambda, \sum_{D\in
\D_{G,X}}\langle\varphi_D,\lambda\rangle D)$. Define the
homomorphism $\chi:\Pic_G(X)\rightarrow \X(T)\times
\Z^{\D_{G,X}}/\X_{G,X}$ by $\chi(\L)=(\lambda, \sum_{D\in
\D_{G,X}}\ord_D(\sigma)D)$, where $\sigma$ is a $B$-semiinvariant
rational section of $\L$ and and $\lambda$ is the weight of
$\sigma$. Extend $\chi$ to a linear map
$\chi:H^2_K(X,\R)\rightarrow \t(\R)\oplus
\R^{\D_{G,X}}/\a_{G,X}(\R)$.  For $v\in H^2_K(X,\R)$ put
$$P(\chi,v)=\chi_0(v)+\{\lambda\in \a_{G,X}(\R)|\langle\lambda,\varphi_D\rangle\geqslant -\chi_D(v),
\forall D\in\D_{G,X}\},$$ where $(\chi_0,\sum_{D\in \D_{G,X}}\chi_D
D)$ is a lifting of  $\chi$. It follows from assertion 2 of
Proposition \ref{Prop:8.6} and Lemma \ref{Lem:8.5} that $2\pi
i\im\psi=P(\chi,[\widetilde{\omega}])$.

\begin{proof}[Proof of Theorem \ref{Thm:8.1}]
Let $X_j$ denote the manifold $M$ with the complex structure
$I_j,j=1,2$. By assertion 3 of Proposition \ref{Prop:8.6},
$\X_{G,X_1}=\X_{G,X_2}$. In the sequel we write $\X$ instead of
$\X_{G,X_j}$ and $\a$ instead of $\a_{G,X_j}$. Put
$\D_j:=\D_{G,X_j}$. Let $\chi^j$ denote the map
$H^2_K(M,\R):\rightarrow (\t(\R)\oplus \R^{\D_j})/\a(\R)$ defined
above, and $\overline{\chi}^j=\chi_0^j+\sum_{D\in \D_{G,X}}\chi_D^i
D:H^2_K(M,\R)\rightarrow \t(\R)\oplus\R^{\D_j}$ be a rational
lifting of $\chi^j$.

{\it Step 1.} Let us reduce the proof to the case when
$[\widetilde{\omega}]$ is rational. By the above,
$$2\pi i\im\psi=P(\chi^i,[\widetilde{\omega}]):
=\chi^i_0([\widetilde{\omega}])+\{\lambda\in \a(\R)|
\langle\varphi_D,\lambda\rangle\geqslant
-\chi_D^i([\widetilde{\omega}]), \forall D\in \D_i\}, i=1,2.$$ Note
that the projections of $\chi_0^1([\widetilde{\omega}]),
\chi^2_0([\widetilde{\omega}])$ to $\t(\R)/\a(\R)$ coincide. Thus,
possibly after  modifying  $\overline{\chi}^1$, we may assume that
$\chi^1_0(v)=\chi^2_0(v)$ for any $v\in V'$, where $V'$ is a
rational subspace $V'\subset H^2_K(M,\R)$ such that
$[\widetilde{\omega}]\in V'(\R)$.

The following lemma implies that there is a sequence
$\widetilde{\omega}_k\in \Pic_G(X)\otimes_\Z \Q$ such that
$P^1(\chi^1,[\widetilde{\omega}_k])=P^2(\chi^2,[\widetilde{\omega_k}])$
and $[\widetilde{\omega}_k]\rightarrow [\widetilde{\omega}]$.

\begin{Lem}\label{Lem:8.7}
Let $V,\a$ be finite dimensional  vector spaces over $\Q$,
$\varphi^1_i,\varphi^2_j\in \a^*, \chi^1_i,\chi^2_j\in V^*$, where
$i=\overline{1,k_1},j=\overline{1,k_2}$. Suppose there is $v\in
V(\R)$ such that the polytopes $P^l(v)\subset\a(\R),l=1,2,$ given by
$$P^l(v):=\{\lambda\in \a(\R)| \langle\varphi^l_i,\lambda\rangle\geqslant -\chi^l_i(v),i=\overline{1,k_i}\}.$$
coincide and have dimension $\dim\a(\R)$. Then for any neighborhood
$O$ of $v$ in $V(\R)$ there is $v'\in O\cap V$ such that
$P^1(v')=P^2(v')$.
\end{Lem}
\begin{proof}[Proof of Lemma \ref{Lem:8.7}]
Let $m$ be the number of facets of $P^l(v)$. We may assume that the
facets are defined by the equations
$\langle\varphi^l_i,\lambda\rangle =\chi^l_i(v),i=\overline{1,m}$.
After rescaling $\varphi^1_i,\chi^1_i,i\leqslant l$, we get
$\varphi^1_i=\varphi^2_i, \chi^1_i=\chi^2_i$ for any $i\leqslant m$.
Let $P(v'), v'\in V(\R),$ denote the polytope given by the
inequalities $\langle \varphi^l_i,\lambda\rangle\geqslant
-\chi_i^l(v'), i=\overline{1,m}$. If the hyperplane $\{\lambda|
\langle\varphi^l_j,\lambda\rangle=\chi_j^l(v') \}$ does not
intersect $P(v')$ for $v'=v$, then the same holds for any $v'$ from
a certain neighborhood of $v$.  Now suppose the hyperplane
$\{\lambda|\langle \varphi^l_j, \lambda\rangle=\chi^l_j(v)\}$ meets
$P^l(v)$ at a face $\Gamma$. Let $I$ denote the subset of
$\{1,\ldots,l\}$ consisting of all $i$ such that
$\langle\varphi^l_i,\cdot\rangle-\chi^l_i(v)$ is a facet of $P^l(v)$
containing $\Gamma$. Then there are positive rational numbers
$a_i,i\in I,$ such that $\varphi^l_j=\sum_{i\in I}a_i\varphi^l_i$.
It follows that $\chi_j^l(v)=\sum_{i\in I}a_i\chi^l_i(v)$. This
equation defines a subspace $V_j^l\subset V$. Let $V_0$ denotes the
intersection of all subspaces $V_j^l$. Since $v\in V_0$, we see that
$V_0$ is nonzero. It follows that $P^l(v')=P(v'), l=1,2,$ for any
$v'$ from a certain neighborhood of $v$ in $V_0$.
\end{proof}
Suppose we have constructed  diffeomorphisms $\varphi_k:M\rightarrow
M$ such that $\varphi_k^*(I_2)=I_1,
\varphi_k([\widetilde{\omega}_k])=[\widetilde{\omega}_k]$. Note that
$\psi_{kl}:=\varphi_k^{-1}\circ\varphi_l$ is a $G$-equivariant
(polynomial) automorphism of $X_1$. By the definition of
$\psi_{kl}$,
\begin{equation}\label{eq:8.1}\psi_{kl}^*([\widetilde{\omega}_l])=\varphi_k^{*-1}([\widetilde{\omega}_l]).\end{equation}
By \cite{Knop8}, the group $\Aut^G(X_1)$ is algebraic. Clearly,
$\Aut^G(X_1)^\circ$ acts trivially on $H^2_K(M,\R)$. So replacing
the sequence $\varphi_k$ with a subsequence, we may assume that the
isomorphism $\psi_{kl}^*$ is the identity on $H^2_K(M,\R)$ for all
$k,l$. Since $[\widetilde{\omega}_k]\rightarrow
[\widetilde{\omega}]$, it follows from (\ref{eq:8.1}) that
$\varphi_k^*([\widetilde{\omega}])=[\widetilde{\omega}]$.

{\it Step 2.} Multiplying $\widetilde{\omega}$ by a sufficiently
large integer $m$, we may assume that
$\widetilde{\omega}=\widetilde{\omega}_{\L^i}$ for a very ample line
$G$-bundle $\L^i$. Let $\widetilde{X}_i$ denote the affine cone over
$X_i$ corresponding to $\L^i$. Set $\widetilde{G}:=G\times
\K^\times$. The variety $\widetilde{X}_i$ has a natural structure of
a (spherical) $\widetilde{G}$-variety. By assertion 1 of Proposition
\ref{Prop:8.6},  $\X^+_{\widetilde{G},\widetilde{X}_i}$ is generated
by integral points in the moment polytope. From \cite{Knop3} it
follows that
$\V_{\widetilde{G},\widetilde{X}_1}=\V_{\widetilde{G},\widetilde{X}_2}$.
Applying Theorem \ref{Thm:1}, we see that
$\widetilde{X}_1\cong^{\widetilde{G}}\widetilde{X}_2$. Therefore
there is a $G$-equivariant isomorphism $\varphi:X_1\rightarrow X_2$
such that $\varphi^*(\L^2)=\L^1$.
\end{proof}

\bigskip

{\Small Massachusetts Institute of Technology, Department of Mathematics, Cambridge MA 02139 USA.

E-mail address: ivanlosev@math.mit.edu}
\end{document}